\documentclass[11pt,reqno]{amsart}
\usepackage{amsmath}
\usepackage{amssymb}
\usepackage{amsfonts}
\usepackage{euscript}
\usepackage[latin1]{inputenc}
\usepackage[T1]{fontenc}
\numberwithin{equation}{section}
\makeatletter\@addtoreset{equation}{section}

\newtheorem {theorem}{Theorem}[section]
\newtheorem {definition}[theorem]{Definition}
\newtheorem {lemma}[theorem]{Lemma}
\newtheorem {proposition}[theorem]{Proposition}
\newtheorem {remark}[theorem]{Remark}

\newtheorem {corollary}[theorem]{Corollary}

\newcommand{\set}[1]{{\left\{{#1}\right\}}}

\newcommand{\abs}[1]{{\left|{#1}\right|}}
\newcommand{\scal}[1]{{\left\langle{#1}\right\rangle}}
\newcommand{\C}{\mathbb C}
\newcommand{\R}{\mathbb R}
\newcommand{\Z}{\mathbb Z}

\newcommand{\K}{\mathbf K}
\newcommand{\Ker}{K}
\newcommand{\Av}{\mathbf{A}}
\newcommand{\La}{\mathbb L}
\newcommand{\fin}{\hfill  
$\square$}

\begin{document}

\title[Landau automorphic functions]{Landau automorphic functions
on $\C^n$ \\ of magnitude $\nu$}

\author{A. Ghanmi $\&$ A. Intissar}
\email{(A.G.) allalghanmi@gmail.com}
  \email{(A.I.)
intissar@fsr.ac.ma}

\thanks{A. G. acknowledges the
financial support of the {\it Arab Regional Fellows Program} within
        a Research Fellow for the full academic year 2006-2007.}

\date{\today}


 \maketitle
\vspace*{-.8cm}
\begin{center}
{\it Department of Mathematics,  Faculty of  Sciences, P.O. Box 1014, \\
   Mohammed V University,  Agdal,  10000 Rabat - Morocco}
   \end{center}
\begin{abstract}
We investigate the spectral theory of the invariant Landau Hamiltonian
\[ 
\La^\nu =  - \frac 12\set{4\sum_{j=1}^n\frac{\partial^2}{\partial z_j \partial \bar z_j} + 2\nu
\sum_{j=1}^n (z_j\frac{\partial}{\partial z_j}-\bar z_j \frac{\partial}{\partial \bar z_j})- \nu^2 |z|^2}
\] 
acting on the space ${\mathcal{F}}^\nu_{\Gamma,\chi}$ of $(\Gamma,\chi)$-automotphic functions on $\C^n$,
constituted of $\mathcal C^\infty$ functions satisfying the functional equation
\[ 
f(z+\gamma) = \chi(\gamma) e^{i\nu\Im m\scal{z,\gamma}}f(z); \quad z\in \C^n,\, \gamma\in\Gamma,
\]
for given real number $\nu>0$, lattice $\Gamma$ of $\C^n$ and a map $\chi:\Gamma\to U(1)$ such that 
the triplet $(\nu,\Gamma,\chi)$ satisfies a Riemann-Dirac quantization type condition.  More precisely,
 we show that  the eigenspace ${\mathcal{E}}^\nu_{\Gamma,\chi}(\lambda)=\set{f\in
{\mathcal{F}}^\nu_{\Gamma,\chi}; \, \, \La^\nu f = \nu(2\lambda+n) f }$; $\lambda\in\C,$  is nontrivial 
if and only if $\lambda=l=0,1,2, \cdots$. In such case, ${\mathcal{E}}^\nu_{\Gamma,\chi}(l)$ is a
finite dimensional vector space whose the dimension is given explicitly by 
\[
\dim{\mathcal{E}}^\nu_{\Gamma,\chi}(l)= \left(\begin{array}{c} n+l-1 \\
l\end{array}\right)({\nu}/{\pi})^n \mbox{vol}(\C^n/\Gamma).
\]
Furthermore, we show that the eigenspace ${\mathcal{E}}^\nu_{\Gamma,\chi}(0)$ associated with the 
lowest Landau level of $\La^\nu$ is isomorphic to the space, ${\mathcal{O}}^\nu_{\Gamma,\chi}(\C^n)$, of 
holomorphic functions on $\C^n$ satisfying 
$$
 g(z+\gamma) = \chi(\gamma) e^{\frac \nu 2 |\gamma|^2+\nu\scal{z,\gamma}}g(z), \eqno{(*)}
$$
that we can realize also as the null space of the differential operator 
$$
\sum\limits_{j=1}\limits^n( \frac{-\partial^2}{\partial z_j\partial \bar z_j} 
+ \nu \bar z_j \frac{\partial}{\partial \bar z_j})
$$ 
acting on  $\mathcal C^\infty$ functions on $\C^n$ satisfying $(*)$.
\end{abstract}

 \section{Notation and statement of main result}

Let $\C^{n}$ be the $n$-complex space endowed with its Hermitian
form $\scal{z,w} = z_1\bar w_1+\cdots +z_n\bar w_n$ and let
$\omega(z,w) = \Im m \scal{z,w}$ be the associated symplectic form.
For given fixed $\nu>0$, we denote by $\La^\nu$ the Landau
Hamiltonian (called also twisted Laplacian \cite{Thangavelu,Ricci}).
It goes back to Landau and describes (for $n=1$) a nonrelativistic
quantum particle moving on the $(x,y)$-plane under the action of an
external constant magnetic field of magnitude $\nu$. The operator
$\La^\nu$ is given explicitly in the complex coordinates $(z_1,z_2,
\cdots , z_n)=z$ by
\begin{equation} \La^\nu =  - \frac 12\set{4\sum_{j=1}^n
\frac{\partial^2}{\partial z_j \partial \bar z_j} + 2\nu
\sum_{j=1}^n (z_j\frac{\partial}{\partial z_j}-\bar z_j
\frac{\partial}{\partial \bar z_j})- \nu^2 |z|^2}\label{LH}
 \end{equation}
and is linked \cite{Thangavelu,Wong} to the sub-Laplacian
\[ {\mathcal{L}} = 4\sum\limits_{j=1}\limits^n \frac{\partial
^2}{\partial z_j
\partial {\bar z_j}}+ 2i \sum\limits_{j=1}\limits^n\Big(z_j\frac{\partial}
{\partial z_j}-{\bar z_j}\frac{\partial}{\partial {\bar
z_j}}\Big)\frac{\partial }{\partial t} + \abs{z}^2\frac{\partial^2
}{\partial t^2}, \]
 on the Heisenberg group
$H^{2n+1}_{eis}=\C^n_z\times\R_t$ through the Fourier transform in
$t$.  Such operator plays an important role in many different
contexts such as Feynman path integral, oscillatory stochastic
integral and theory of lattices electrons in uniform magnetic field
(see Bellissard \cite{Bellissard} and references therein). For what
is the spectral properties, it is known that $\La^\nu$ is a
selfadjoint elliptic differential operator on $L^2(\C^{n};dm)$, the
usual Hilbert space of square integrable functions on $\C^n$ with
respect to the Lebesgue measure $dm$. Its spectrum is purely
discrete and given by the eigenvalues (Landau levels)
\begin{equation}
\nu(2l+n), \qquad l=0,1,2, \cdots ,
\end{equation}
which occur with infinite multiplicities
\cite{AvronHerbestSimon,Novikov,Shigekawa}. Additonal spectral
properties relevant for our purpose are recalled in Section 2.

In this paper, we consider the action of the operator $\La^\nu$ on
some appropriate functional spaces. Let $\Gamma$ be a full rank
lattice of $\C^n$ (i.e., a
 discrete subgroup of rank $2n$ of the additive group $\R^{2n}=\C^n$) so that $\C^n/\Gamma$
 is compact, and  $\chi$ be a given map
 \[ \chi : \Gamma \longrightarrow U(1)=\set{\lambda\in \C; ~ |\lambda|=1}.\]
To the given data ($\nu,\Gamma,\chi$), we then associate the space
${\mathcal{F}}^\nu_{\Gamma,\chi}$ of
 $\mathcal{C}^\infty$ functions $f$ on $\C^n$ that satisfy the functional equation
 \begin{equation}
f(z+\gamma) = \chi(\gamma) e^{i\nu\omega(z,\gamma)}f(z)
\label{ChiGammaPeriodic}
 \end{equation} for all $z\in\C^n$ and $ \gamma\in \Gamma$, as well as the space
 ${\mathcal{O}}^\nu_{\Gamma,\chi}(\C^n)$
 of holomorphic functions $g$ on $\C^n$, $g\in \mathcal{O}(\C^n)$, satisfying the following
 functional equation
  \begin{equation}
g(z+\gamma) = \chi(\gamma) e^{\frac \nu 2
|\gamma|^2+\nu\scal{z,\gamma}}g(z) 
\label{ChiGammaPeriodicEntire}
 \end{equation}
for all $z\in \C^n$ and $ \gamma\in \Gamma$. Then, it will be shown
that the space ${\mathcal{F}}^\nu_{\Gamma,\chi}$ (or also
${\mathcal{O}}^\nu_{\Gamma,\chi}(\C^n)$) is a nonzero complex
vector space if and only if the
  triplet $(\nu,\Gamma, \chi)$ satisfies the following
$(RDQ)$ condition $$\chi(\gamma_1+\gamma_2)=
\chi(\gamma_1)\chi(\gamma_2)e^{i\nu\omega(\gamma_1,\gamma_2)}
\eqno(RDQ) $$ for  every $\gamma_1,\gamma_2\in \Gamma$ (see Proposition
\ref{NoTrivial1}). In this case and owing to the fact that the
Landau Hamiltonian $\La^\nu$ leaves invariant
 the space ${\mathcal{F}}^\nu_{\Gamma,\chi}$ (see Proposition
\ref{Proposition31}), we can consider its restriction to
${\mathcal{F}}^\nu_{\Gamma,\chi}$ that we shall denote by
$\La^\nu_{\Gamma,\chi}$ and therefore consider the associated
eigenvalue problem
$  \La^\nu_{\Gamma,\chi} f = \nu(2\lambda+n)  f 
$ in ${\mathcal{F}}^\nu_{\Gamma,\chi}$ with $\lambda \in \C$. Hence
by ${\mathcal{E}}^\nu_{\Gamma,\chi}(\lambda)$ let denote the
corresponding eigenspace, i.e.,
 \begin{equation}
 {\mathcal{E}}^\nu_{\Gamma,\chi}(\lambda)=
 \set{f \in {\mathcal{F}}^\nu_{\Gamma,\chi} ;
 \quad \La^\nu_{\Gamma,\chi} f = \nu(2\lambda+n) f}.
 \label{EigenPeriodic}
 \end{equation}
and let make for instance the following
 \begin{definition}
\label{definition} Assume $(RDQ)$ to be satisfied by the triplet
$(\nu,\Gamma, \chi)$.
\begin{itemize}
 \item[i)] We call
 ${\mathcal{F}}^\nu_{\Gamma,\chi}$  the space of
 $(\Gamma,\chi)$-automorphic functions  on $\C^n$ of magnitude $\nu$.
 \item[ii)] We call ${\mathcal{E}}^\nu_{\Gamma,\chi}(\lambda)$ the space of
 Landau $(\Gamma,\chi)$-automorphic functions of magnitude
 $\nu$ at the level $\lambda$. The particular one  \begin{equation}
 {\mathcal{E}}^\nu_{\Gamma,\chi}(0):=
 \set{f \in {\mathcal{F}}^\nu_{\Gamma,\chi} ;
 \quad \La^\nu_{\Gamma,\chi} f = n\nu  f}.
 \label{EigenPeriodic0}
 \end{equation} corresponding to
 $\lambda=0$ is called the "fundamental space of $(\Gamma,\chi)$-ground
states".
 \item[iii)] We call ${\mathcal{O}}^\nu_{\Gamma,\chi}(\C^n)$ the space of
 holomorphic $(\Gamma,\chi)$-automorphic functions of magnitude
 $\nu$ or also  $(\Gamma,\chi)$-theta functions on $\C^n$.
 \end{itemize}
\end{definition}

 The objective of the present paper
is to investigate the spectral analysis of the involved eigenspaces
 ${\mathcal{E}}^\nu_{\Gamma,\chi}(\lambda)$.
 Namely, the main result to which is aimed this paper is the
 following

 \vspace*{.25cm}
 \noindent {\bf Main Theorem.} {\it Assume the $(RDQ)$ condition to be satisfied
by the triplet $(\nu,\Gamma, \chi)$ and let
${\mathcal{E}}^\nu_{\Gamma,\chi}(\lambda)$ and
${\mathcal{O}}^\nu_{\Gamma,\chi}(\C^n)$ be the functional spaces
defined above.  Then
\begin{itemize}
    \item[i)] The eigenspace ${\mathcal{E}}^\nu_{\Gamma,\chi}(\lambda)$ is a
nonzero vector space if and only if  $\lambda$ is a positive
integer $\lambda=l=0,1,2, \cdots$.
    \item[ii)] For every fixed positive integer $l=0,1,2, \cdots$, the space
 \[{\mathcal{E}}^\nu_{\Gamma,\chi}(l)=\set{f; \, f\in {\mathcal{F}}^\nu_{\Gamma,\chi},
  \,\, \La^\nu_{\Gamma,\chi}f=\nu(2l+n)f}\] is a finite dimensional vector space whose the
  dimension is given explicitly by the formula
\begin{equation} \label{DimFormulaIntr} \dim{\mathcal{E}}^\nu_{\Gamma,\chi}(l)=
\frac{\Gamma(n+l)}{\Gamma(n)l!}({\nu}/{\pi})^n
\mbox{vol}(\C^n/\Gamma) ,
\end{equation}
 where $\Gamma(x)$ is the usual gamma function and $\mbox{vol}(\C^n/\Gamma)$ denotes the Lebesgue volume
of a fundamental domain of the lattice $\Gamma$.
    \item[iii)] The fundamental space of $(\Gamma,\chi)$-ground states,
 ${\mathcal{E}}^\nu_{\Gamma,\chi}(0)$, is isomorphic to the space
${\mathcal{O}}^\nu_{\Gamma,\chi}(\C^n)$ with $f \longmapsto
g=e^{\frac\nu 2|z|^2}f$ from ${\mathcal{E}}^\nu_{\Gamma,\chi}(0)$
onto ${\mathcal{O}}^\nu_{\Gamma,\chi}(\C^n)$ as isomorphism map.
\end{itemize}
}

\vspace*{.15cm}

For the establishment of our main result, we have make use of the
explicit description of the spectral analysis of the operator
$\La^\nu$. The computation of the dimension of the eigenspaces $
{\mathcal{E}}^\nu_{\Gamma,\chi}(l) $ is done \`a la Selberg
\cite{Selberg,Hejhal}. In fact, we determinate the traces of
integral operators associated with $(\Gamma,\chi)$-automorphic kernel
functions obtained by averaging reproducing kernels of the free
$L^2$-eigenspaces of $\La^\nu$.

\begin{remark} \label{RemStability} \quad
\begin{itemize}
    \item[a)] Owing to Proposition \ref{Proposition31}, the statement i) in the
main theorem shows that the operator $\La^\nu$ and its restriction
$\La^\nu_{\Gamma,\chi}
:={\La^\nu}|_{{\mathcal{F}}^\nu_{\Gamma,\chi}}$  have the same
spectrum (i.e., stability of the spectrum under perturbation by the
lattice $\Gamma$). However, the degeneracy of the eigenvalues
becomes finite.

\item[b)] Note that the dimension of the space of
 Landau $(\Gamma,\chi)$-automorphic functions ${\mathcal{E}}^\nu_{\Gamma,\chi}(l)$ is
independent of the multiplier $\chi$. It can also be noted that all
the eigenspaces have the same dimension when $n=1$. While for $n\geq
2$ the dimension of the spaces ${\mathcal{E}}^\nu_{\Gamma,\chi}(l)$
growths polynomially in $l$. Namely, we have
\[ \dim{\mathcal{E}}^\nu_{\Gamma,\chi}(l) \sim C l^{n-1} \quad
\mbox{as } l\to +\infty \] for certain constant $C>0$.
\end{itemize}
\end{remark}

The outline of the paper is as follows. In Section 2, we collect and
review some needed background on the spectral theory of the Landau
Hamiltonian $\La^\nu$ acting on the free Hilbert space
$L^2(\C^n;dm)$. Section 3 is devoted to give some basic properties
of the spaces ${\mathcal{F}}^\nu_{\Gamma,\chi}$ and
${\mathcal{O}}^\nu_{\Gamma,\chi}(\C^n)$. Mainely, we prove the
equivalence of the nontriviality of such spaces to the $(RDQ)$
condition and moreover we
 explicit the expression of the reproducing kernel of
${\mathcal{O}}^\nu_{\Gamma,\chi}(\C^n)$ as well as its dimension. In
Section 4, we investigate some general properties of a class of
$(\Gamma,\chi)$-automorphic kernel functions on $\C^n\times\C^n$ of
magnitude $\nu$ that are essential for our purpose.
 In Section 5, we present the proof of our main result. The latest section deals with some concluding remarks.

\vspace*{.15cm} We conclude this introduction by providing an
example of triplet  $(\nu,\Gamma, \chi)$ satisfying the $(RDQ)$
condition.  For $\Gamma$ being a lattice in $\C=\R^2$, we denote by
$S_{_\Gamma}$ its cell area and we set
$\nu_{_\Gamma}=\pi/S_{_\Gamma}$. Let $\chi_{_\Gamma}$ be the
Weierstrass pseudo-character defined on $\Gamma$ by
$\chi_{_\Gamma}(\gamma)=+1$ if $\gamma /2 \in \Gamma$ and
$\chi_{_\Gamma}(\gamma)=-1$ otherwise \cite[page
103]{Jones-Singerman}. Then, it can be shown that
$(\nu_{_\Gamma},\Gamma, \chi_{_\Gamma})$ satisfies $(RDQ)$ and that
we have
$\dim{\mathcal{E}}^{\nu_{_\Gamma}}_{\Gamma,\chi_{_\Gamma}}(l)=1$ for
every $l=0,1,2, \cdots$. Moreover, one can build generator of each
${\mathcal{E}}^{\nu_{_\Gamma}}_{\Gamma,\chi_{_\Gamma}}(l)$ involving
basically the modified Weierstrass sigma function.

 \section{Background on spectral theory of the operator $\La^\nu$}

 \noindent We begin with an invariance property of the Landau Hamiltonian
 $\La^\nu$. For this, let $G=U(n) \rtimes \C^n$ be the
solvable semidirect product of the unitary group $U(n)$ with the
additive group $(\C^n,+)$. Such group is also realized as
\[ G=\set{g=\left(
\begin {array} {c c} a & b \\ 0 & 1
\end {array} \right ); \qquad a\in U(n) , ~ b\in \C^n}.\]
and acts transitively on $\C^n$ by the holomorphic mappings $z
\longmapsto g.z :=az+b$, which can be extended to $L^2(\C^n;dm)$ by
considering
\begin{equation}
[T^\nu_{g}f](z):=j_\nu(g,z) f(g.z), \label{MT01}
 \end{equation} where the involved factor $j_\nu(g,z)$ is given by
\begin{equation}j_\nu(g,z) = e^{i\nu\omega(z,g^{-1}.0)}.\label{MT02}
 \end{equation}
 We then assert
\begin{proposition}\label{Proposition31} \quad

\begin{itemize}
    \item[i)] For every $g_1,g_2\in G$, we have the chain rule
 \begin{equation}\label{ChainRule}
 j_{\nu}(g_1g_2,z) =
 e^{i\nu\omega(g^{-1}_1.0,g_2\cdot 0)}  j_{\nu}(g_1,g_2.z)
 j_{\nu}(g_2,z).
\end{equation}
  \item[ii)] The transformation $T^\nu$ defines a projective representation of the group $G$ on
 the Hilbert space $L^2(\C^n;dm)$. That is
  \begin{itemize}
    \item[a)] It is a unitary transformation on $L^2(\C^n;dm)$ for every $g\in G$.
    \item[b)] For all $g_1, g_2\in G$, we have $T^\nu_{g_1 g_2} =
     e^{i\phi_\nu(g_1,g_2)}T^\nu_{g_2}\circ T^\nu_{g_1}$, where the phase factor is given here by $\phi_\nu(g_1,g_2)=
     \nu\omega(g^{-1}_1.0,g_2\cdot 0)$.
   \item[c)] The map $g \mapsto T^\nu_{g}f$ from $G$ into $L^2(\C^n;dm)$ is a
continuous map for every fixed $f\in L^2(\C^n;dm)$. \end{itemize}
 \item[iii)] The Landau Hamiltonian $\La^\nu$ is $T^\nu$-invariant in the sense
 that for every $g\in G$ we have
  $T^\nu_{g} \La^\nu = \La^\nu T^\nu_{g}$.
  \end{itemize}
\end{proposition}

\noindent The proof of such proposition can be handled by
straightforward computation. For iii), one can also refer to
\cite{GI} for an intrinsic different approach.

\vspace*{.2cm} \noindent Additional needed spectral properties of
$\La^\nu $ are summarized in the following

\begin{proposition} \label{Proposition32}\quad

\begin{itemize}
    \item[i)] For fixed $\lambda\in \C$, the space of  radial functions $f$
   solution  of \,  $\La^\nu f =\nu(2\lambda+n) f $ is one
   dimensional, and it is generated  by
\begin{equation}\label{radialsolution}
\varphi_\lambda(z) = e^{-\frac\nu 2|z|^2} {_1F_1}(-\lambda; n;
\nu|z|^2),
\end{equation}
  where ${_1F_1} (a;c;x)=1+\frac{a}{c}\frac{x}{1!}+\frac{a(a+1)}{c(c+1)}\frac{x^2}{2!}+\cdots $
  is the usual confluent hypergeometric function.
\item[ii)] The function $\varphi_\lambda$ given by (\ref{radialsolution}) is
bounded if and only if $\lambda$ is a
    positive integer $l$; $l=0,1,2, \cdots $.
\item[iii)] The spectrum of ~$\La^\nu$ acting on $L^2(\C^n;dm)$ is
 discrete and given by the so-called Landau levels
$\nu(2l+n)$; $l=0,1,2, \cdots , $
     where each eigenvalue occurs with infinite multiplicity. Furthermore, we
     have the following  orthogonal decomposition in Hilbertian subspaces
\begin{equation} \label{freeDecomp}
 L^2(\mathbb{C}^{n};dm) = \bigoplus_{l=0}^\infty \mathcal{A}^{2,\nu}_{l}(\C^n)
\end{equation}
    where $\mathcal{A}^{2,\nu}_{l}(\C^n)=\set{ f; ~ f\in L^2(\mathbb{C}^{n};dm) ~ \mbox{ and } ~  \La^\nu f = \nu(2l+n) f}.$
\item[iv)] The $L^2$-eigenprojector kernel
of the $L^2$-eigenspace $\mathcal{A}^{2,\nu}_{l}(\C^n)$ is given
explicitly by
    the following closed formula
\begin{equation}\label{RepK}
{\mathcal{K}}^\nu_l(z,w) = e^{i\nu\omega(z,w)}Q^\nu_l(|z-w|) ,
\end{equation}
    where we have set
\[\label{Exp22}
 Q^\nu_l(|z-w|) =
\frac{\Gamma(n+l)}{\Gamma(n)l!} ({\nu}/{\pi})^n e^{-\frac\nu
2|z-w|^2} {_1F_1}(-l; n; \nu |z -w|^2),
\]
and then  satisfies the following "$G$-invariance property"
\begin{equation}\label{InvProp}
{\mathcal{K}}^\nu_l(z,w)=
e^{i\nu\omega(z,g^{-1}.0)}{\mathcal{K}}^\nu_l(g.z,g.w)e^{-i\nu\omega(w,g^{-1}.0)}
\end{equation}
    for every $g\in G$ and $z,w\in \C^n$.
\end{itemize}
\end{proposition}

\noindent{\it Proof.} To check the statement i), we write $\La^\nu$
in polar coordinates $z=r\theta$; $r\geq 0$ and $\theta$ in the
$(2n-1)$-dimensional unit sphere $S^{2n-1}$,
\[-2\La^\nu = \frac{\partial^2}{{\partial r}^2 }+ \frac{2n-1}{r}\frac{\partial}{\partial r}
+\frac 1{r^2}\Delta_{S^{2n-1}} +2\nu(L_\theta-\overline{L_\theta})
-\nu^2r^2 ,\] where $\Delta_{S^{2n-1}}$ denotes the Laplace-Beltrami
operator on $S^{2n-1}$ and $L_\theta$ is the tangential component of
the complex Euler operator
$E=\sum\limits_{j=1}\limits^nz_j\frac{\partial}{\partial
z_j}=\frac{r}{2}\frac{\partial}{\partial r}+L_\theta$. So that the
differential equation $\La^\nu f =\nu(2\lambda+n) f $, for radial
solutions $\phi(r)$, reduces to the following
\[ \frac{d^2\phi}{{d r}^2 }+ \frac{2n-1}{r}\frac{d\phi}{d r} -\nu^2r^2\phi +2 \nu(2\lambda+n)\phi= 0 . \]
 Next, making use of the appropriate
change of function $\phi(r)=e^{- x/2}y(x) $ with $x=\nu r^2$, we see
that the previous equation leads to the following ordinary
differential equation \cite[page 193]{NO}
\[
xy^{''}+(n-x)y^{'}+\lambda y = 0
\]
whose regular solution at $x=0$ is the confluent hypergeometric
function  ${_1F_1}(-\lambda;n;x)$.

The assertion ii) is clear for $\lambda=0$. For $\lambda\ne 0$, we
use the asymptotic behavior of the confluent hypergeometric function
given by \cite[page 332]{NO}
\[
{_1F_1}(a; c; x) = \Gamma(c)\set{\frac{(-x)^{-a}}{\Gamma(c-a)} +
\frac{e^{x} x^{a-c}}{\Gamma(a)}  } \Big(1+O(\frac 1 x )\Big) ~~
\]
as $x \rightarrow +\infty$. Hence for $a=-\lambda$, $c=n$ and
$x=\nu|z|$ in above, we obtain \[
 \lim\limits_{|z| \rightarrow +\infty} e^{- \frac\nu
2 |z|^2} {_1F_1}(-\lambda; n; \nu |z|^2) = \lim\limits_{|z|
\rightarrow +\infty} \frac{\Gamma(n)}{\Gamma(-\lambda)}
(\nu|z|^2)^{-(n+\lambda)}e^{\nu |z|^2}  .
\]
From which, we conclude that $\varphi_\lambda$ is bounded if and
only if $\lambda=1,2, \cdots$.

The result in iii)  is  well known and the reader can refer for
example to \cite{AvronHerbestSimon}. While the proof of iv) is
contained in \cite{GI} and can be handled in a similar way as in
\cite{AskourIntissarMouayn}.
 \fin\\

 In the next section, we study some basic properties of the space of
 $(\Gamma,\chi)$-automorphic functions that the Landau
 Hamiltonian $\La^\nu$ will act on.

\section{Basic properties of the space
      ${\mathcal{F}}^\nu_{\Gamma,\chi}$ and associated spaces}

 Recall that for given data
$\nu > 0$, $\Gamma $ a lattice of $\R^{2n}=\C^n$ of rank $2n$ and
$\chi$ a mapping from $\Gamma$ to the unit circle $\set{\lambda \in
\C; ~ |\lambda|=1}= U(1)$, we have associated the functional space
\begin{equation}
 {\mathcal{F}}^\nu_{\Gamma,\chi}=\set{f \in {\mathcal{C}}^\infty(\C^n) ;
 \quad f(z+\gamma) = \chi(\gamma) e^{i\nu\omega(z,\gamma)}f(z)  } .
 \label{AutForms}
 \end{equation}
The following proposition gives sufficient and necessary condition
on the triplet $(\nu,\Gamma,\chi)$ in order that
${\mathcal{F}}^\nu_{\Gamma,\chi}$ is a nonzero space. Namely, we
have

\begin{proposition}\label{NoTrivial1}
The complex vector space ${\mathcal{F}}^\nu_{\Gamma,\chi}$ is a nonzero space if and only if the triplet $(\nu,\Gamma, \chi)$ satisfies
the following condition $$\chi(\gamma_1+\gamma_2)=
\chi(\gamma_1)\chi(\gamma_2)e^{i\nu\omega(\gamma_1,\gamma_2)}
\eqno(RDQ) $$ for  every $\gamma_1,\gamma_2\in \Gamma$.
 In this case ${\mathcal{F}}^\nu_{\Gamma,\chi}$ is an
 infinite dimensional complex vector space.
 \end{proposition}

  \begin{remark} \label{RemarkRDQC1}
Under the $(RDQ)$ condition, the map $\chi$ satisfies the following
properties:
\begin{equation} \chi(0) = 1 \quad\mbox{and}\quad  \quad
\chi(-\gamma) =\overline{\chi(\gamma)} . \label{ChiProperties}
\end{equation} Also, by interchanging the roles of $\gamma_1$ and
$\gamma_2$ in $(RDQ)$ and using  the fact that the symplectic form
$\omega(\cdot,\cdot)$ is antisymmetric, we have necessarily
 \begin{equation}
{\nu}\omega(\gamma_1,\gamma_2) \in {\pi} \Z
    \label{ARDQC} \end{equation}
for every $\gamma_1,\gamma_2\in \Gamma$.
  \end{remark}

 \begin{remark} \label{RemarkRDQC2}
  The $(RDQ)$ condition is equivalent to that the complex valued
  function $J_{\nu,\chi}$ defined on $\Gamma\times\C^n$ by
   $J_{\nu,\chi}(\gamma,z):=\chi(\gamma )e^{i\nu \omega(z,\gamma)}$ is an
   automorphy factor satisfying the cocycle identity,
   \[J_{\nu,\chi}(\gamma_1+\gamma_2,z)=J_{\nu,\chi}(\gamma_1,z+\gamma_2)J_{\nu,\chi}(\gamma_2,z).\] Therefore
   $\phi_\gamma(z;v) := (z+\gamma; \chi(\gamma )e^{i\nu
   \omega(z,\gamma)}.
   v)$ defines an action of $\Gamma$ on $\C^n\times\C$ and the
   associated quotient space $(\C^n\times\C)/\Gamma$ is a line
   bundle over the torus $\C^n/\Gamma$ with fiber $\C=\tau^{-1}([z])$, where
   the  projection map
   $\tau: (\C^n\times\C)/\Gamma\longrightarrow \C^n/\Gamma$ is the natural one induced from the canonical projection
   $ \pi: \C^n\longrightarrow \C^n/\Gamma$.
   Thus, one can regard the space  ${\mathcal{F}}^\nu_{\Gamma,\chi}$ as the
   space of $\mathcal{C}^\infty$ sections of the above line bundle
   over the complex torus $\C^n/\Gamma$
   and  therefore   ${\mathcal{F}}^\nu_{\Gamma,\chi}$ is of
   infinite dimension.
  \end{remark}

      \begin{remark} \label{RemarkRDQC11}
The abbreviation  $(RDQ)$ is used to refer to "Riemann-Dirac
Quantization" condition. Indeed, the pair $(H,E)$, with
$H(z,w):=(\nu/\pi)\scal{z,w}$ and $E(z,w):=\Im m H(z,w)$, satisfies
the Riemann condition \cite{GriffithsHarris,Satake78}, 
 and then the considered complex torus is an abelian
variety.  Also, the condition (\ref{ARDQC}) is that called in
Quantum Mechanics Dirac quantization.
  \end{remark}

\noindent{\it Proof of Proposition \ref{NoTrivial1}.} The proof of
"{only if}" follows by assuming that
${\mathcal{F}}^\nu_{\Gamma,\chi}$ is nontrivial space and next by
computing $f(z+\gamma_1+\gamma_2)$ in two manners, for a given nonzero function
 $f\in
{\mathcal{F}}^\nu_{\Gamma,\chi}$. Indeed, we have
\begin{equation} f(z+\gamma_1+\gamma_2) \stackrel{}{=}
\chi(\gamma_1+\gamma_2)e^{i\nu\omega(z,\gamma_1+\gamma_2)} f(z)
\label{RHS1}
\end{equation}
and also
\begin{align} f(z+\gamma_1+\gamma_2)&=  f([z+\gamma_1]+\gamma_2)=
\chi(\gamma_2)e^{i\nu\omega(z+\gamma_1,\gamma_2)}f(z+\gamma_1)
\nonumber
\\
&=  \chi(\gamma_1)\chi(\gamma_2)e^{i\nu\omega(\gamma_1,\gamma_2)}
e^{i\nu\omega(z,\gamma_1+\gamma_2)}f(z).  \label{RHS2}
\end{align}
Next, by equating the right hand sides of (\ref{RHS1}) and
(\ref{RHS2}) and using the fact that
 $f$ is not identically zero, we conclude that
\[\chi(\gamma_1+\gamma_2) =
\chi(\gamma_1)\chi(\gamma_2)e^{i\nu\omega(\gamma_1,\gamma_2)}
\eqno(RDQ)\] for all $ \gamma_1,\gamma_2\in \Gamma$.

  For the converse, we may use Remark \ref{RemarkRDQC2}.
  But for readers whom are not familiar to such language a direct proof can be given.
  In fact by classical analysis, one can pick any arbitrary non
  zero $\mathcal{C}^\infty$ function $\psi$ on $\C^n$  such that
$\mbox{Supp}\psi\subset \Lambda(\Gamma)$, where $\Lambda(\Gamma)$ is
a fundamental domain of $\Gamma$ in $\C^n$. Then we can use the
following lemma to show that the space
${\mathcal{F}}^\nu_{\Gamma,\chi}$ is a nonzero space and it is of
infinite dimension.
\begin{lemma} \label{GammaChiPerio} Suppose that the  condition
  $(RDQ)$ holds. Let $\psi$ be a compactly supported $\mathcal{C}^\infty$ function such that
  $Supp \psi\subset \Lambda(\Gamma)$, and denote by ${\mathcal{P}}^\nu_{\Gamma,\chi}\psi$ the
 $(\Gamma,\chi)$-periodization (\`a la Poincar\'e) of $\psi$ given
by
\[
[{\mathcal{P}}^\nu_{\Gamma,\chi}\psi](z)=\sum\limits_{\gamma\in
\Gamma}
\overline{\chi(\gamma)}e^{-i\nu\omega(z,\gamma)}\psi(z+\gamma)
\]
or equivalently (in view of Remark \ref{RemarkRDQC1}) by
\[
[{\mathcal{P}}^\nu_{\Gamma,\chi}\psi](z)=\sum\limits_{\gamma\in
\Gamma} \chi(\gamma)e^{i\nu\omega(z,\gamma)}\psi(z-\gamma).
\]
Then, we have

i) The function ${\mathcal{P}}^\nu_{\Gamma,\chi}\psi$ is a nonzero
$\mathcal{C}^\infty$ function on $\C^n$.

ii) The function ${\mathcal{P}}^\nu_{\Gamma,\chi}\psi$ belongs to
the space ${\mathcal{F}}^\nu_{\Gamma,\chi}$. \fin
\end{lemma}

\noindent{\it Proof of Lemma \ref{GammaChiPerio}.} i) By
construction the Poincar\'e series
$[{\mathcal{P}}^\nu_{\Gamma,\chi}\psi](z)$ is well defined as
$\mathcal C^\infty$ function on $\C^n$. Furthermore, since $\Gamma$
is discrete, $Supp \psi\subset \Lambda(\Gamma)$ and $\chi(0)=1$,  it
is easy to see that for every $z\in Supp \psi$, we have
\[ [{\mathcal{P}}^\nu_{\Gamma,\chi}\psi](z)=
\psi(z).\]

ii) For every $\gamma \in \Gamma$ and $z\in \C^n$, we have
\begin{align*}
[{\mathcal{P}}^\nu_{\Gamma,\chi}\psi](z+\gamma)&= \sum\limits_{h\in
\Gamma} \chi(h)e^{i\nu\omega(z+\gamma ,h)}\psi([z+\gamma]-h)
\\
&\stackrel{h=\gamma+k}{=} \sum\limits_{k\in \Gamma}
\chi(\gamma+k)e^{i\nu\omega(z+\gamma ,\gamma+k)}\psi(z-k) .
\end{align*}
Therefore, using the $(RDQ)$  condition $\chi(\gamma+k) =
\chi(\gamma)\chi(k)e^{i\nu\omega(\gamma,k)}$, we get
\[
[{\mathcal{P}}^\nu_{\Gamma,\chi}\psi](z+\gamma) =
 \chi(\gamma)e^{i\nu\omega(z,\gamma)} \sum\limits_{k\in \Gamma}
\chi(k)e^{2i\nu\omega(\gamma,k)}e^{i\nu\omega(z,k)}\psi(z-k)
\]
Next, since  $\nu\omega(\gamma,k)\in \pi\Z$ for all
$\gamma,k\in\Gamma$, it follows that $e^{2i\nu\omega(\gamma,k)}=1$
and thus \begin{align*}
[{\mathcal{P}}^\nu_{\Gamma,\chi}\psi](z+\gamma) &=
\chi(\gamma)e^{i\nu\omega(z,\gamma)}\sum\limits_{k\in \Gamma}
\chi(k)e^{i\nu\omega(z,k)}\psi(z-k)
\\&=
\chi(\gamma)e^{i\nu\omega(z,\gamma)}
[{\mathcal{P}}^\nu_{\Gamma,\chi}\psi](z).
\end{align*}
The proof of Lemma \ref{GammaChiPerio} is finished.\fin

\vspace*{.2cm}
 Now noting that for every $f_1, f_2 \in {\mathcal{F}}^\nu_{\Gamma,\chi}$ the product function
 $f_1(z) \overline{f_2}(z) $ is a $\Gamma$-periodic on
 $\C^n$, then $\int_{\C^n/\Gamma}f_1(z) \overline{f_2}(z)
dm(z)$ makes sense, and therefore we can equip
${\mathcal{F}}^\nu_{\Gamma,\chi}$
 with the inner scalar product, to wit
\begin{equation}\label{HSPGamma1}
\scal{f_1, f_2}_\Gamma :=\int_{\C^n/\Gamma}f_1(z) \overline{f_2}(z)
dm(z)= \int_{\Lambda(\Gamma)}f_1(z) \overline{f_2}(z) dm(z),
\end{equation}
 where $\Lambda(\Gamma)$ is any given fundamental domain of the lattice
$\Gamma$. We denote by $L^{2,\nu}_{\Gamma,\chi}$ the completion of
 ${\mathcal{F}}^\nu_{\Gamma,\chi}$ with respect to the norm
 \[| f |_\Gamma =  \sqrt{\scal{f, f}_\Gamma}; \quad f  \in
 {\mathcal{F}}^\nu_{\Gamma,\chi}.\]

\begin{remark} \label{RealizationHil} The Hilbert space $L^{2,\nu}_{\Gamma,\chi}$ can be characterized as the space of all
measurable functions on $\C^n$ that are square integrable on
$\Lambda(\Gamma)$ with respect to the Lebesgue measure $dm$ and
satisfying the functional equation
\[f(z+\gamma) = \chi(\gamma) e^{i\nu\omega(z,\gamma)}f(z)\] for
almost every $z\in \C^n$ and every $\gamma\in \Gamma$.
\end{remark}

 In parallel to ${\mathcal{F}}^\nu_{\Gamma,\chi}$, we can consider
the functional space
\begin{equation}\label{GSpace}
 {\mathcal{G}}^\nu_{\Gamma,\chi}=\set{g \in {\mathcal{C}}^\infty(\C^n) ;
 \quad g(z+\gamma) = \chi(\gamma) e^{\frac \nu 2
|\gamma|^2+\nu\scal{z,\gamma}}g(z)  }.
  \end{equation}
endowed with the norm $||\cdot||_\Gamma$ associated with the Hermitian
scalar product
\begin{equation}\label{HSPGamma2}\scal{\scal{g_1,
g_2}}_\Gamma=\int_{\C^n/\Gamma}g_1(z) \overline{g_2}(z) e^{-\nu
  |z|^2}dm(z)  ,
  \end{equation}
  Then,   we verify that
$({\mathcal{F}}^\nu_{\Gamma,\chi},|\cdot|_\Gamma)$ and
$({\mathcal{G}}^\nu_{\Gamma,\chi},||\cdot||_\Gamma)$ are isometric
pre-Hilbertian spaces through the mapping
\begin{equation}\label{TransformationG} f\in {\mathcal{F}}^\nu_{\Gamma,\chi} \longmapsto \mathfrak{G}
f\in {\mathcal{G}}^\nu_{\Gamma,\chi}; \quad [\mathfrak{G} f](z) =
e^{\frac \nu 2|z|^2} f(z).
\end{equation}
Thus the Landau Hamiltonian $\La^\nu_{\Gamma,\chi}$ acting on
${\mathcal{F}}^\nu_{\Gamma,\chi}$ gives rise to the following second
order differential operator $\Delta^\nu_{\Gamma,\chi}$ acting on
${\mathcal{G}}^\nu_{\Gamma,\chi}$ by means of
 \begin{equation}
 \label{RelationLaplacians} \Delta^\nu_{\Gamma,\chi} g= \frac 12 e^{\frac \nu 2|z|^2}
 \Big[ \La^\nu_{\Gamma,\chi} - n\nu \Big] (e^{-\frac \nu 2|z|^2} g)
 \end{equation}
 for every $g\in {\mathcal{G}}^\nu_{\Gamma,\chi}$.  More precisely, we have
 \begin{equation}\label{GLaplacian}
 \Delta^\nu_{\Gamma,\chi}=
  \sum_{j=1}^n \Big(
\frac{-\partial^2}{\partial z_j \partial \bar z_j}+ \nu
 \bar z_j \frac{\partial}{\partial
\bar z_j} \Big).  \end{equation} Therefore,  describing the spectral
analysis of $\La^\nu_{\Gamma,\chi}$ on
${\mathcal{F}}^\nu_{\Gamma,\chi}$ is equivalent to do it for
$\Delta^\nu_{\Gamma,\chi}$ on the functional space
${\mathcal{G}}^\nu_{\Gamma,\chi}$.

Below, we will focus on the natural subspace $
{\mathcal{O}}^\nu_{\Gamma,\chi}(\C^n)$ of
 ${\mathcal{G}}^\nu_{\Gamma,\chi}$, consisting of
holomorphic functions on $\C^n$, $g\in
  {\mathcal{O}}(\C^n)$,
  satisfying the functional equation
  \[ g(z+\gamma) = \chi(\gamma) e^{\frac \nu 2
|\gamma|^2+\nu\scal{z,\gamma}}g(z)\] for every $z\in \C^n$ and every
$\gamma\in \Gamma$, i.e.,
\begin{equation} {\mathcal{O}}^\nu_{\Gamma,\chi}(\C^n)=
                     \set{g  \in \mathcal{O}(\C^n) ;
                  \quad g(z+\gamma) = \chi(\gamma) e^{\frac \nu 2
                   |\gamma|^2+\nu\scal{z,\gamma}}g(z)}.\label{HolAutForms}
 \end{equation}
Then, according to the isometry ${\mathcal{G}}^\nu_{\Gamma,\chi}
\cong {\mathcal{F}}^\nu_{\Gamma,\chi}$, one concludes easily from
Proposition \ref{NoTrivial1} that $(RDQ)$ is a necessarily condition
to ${\mathcal{O}}^\nu_{\Gamma,\chi}(\C^n)$ be nontrivial subspace
of ${\mathcal{G}}^\nu_{\Gamma,\chi}$. This can be handled directly
as in Proposition \ref{NoTrivial1}. But for the converse, i.e.,
$(RDQ)$ implies ${\mathcal{O}}^\nu_{\Gamma,\chi}(\C^n)\ne \{0\}$,
one has to proceed differently since we do not dispose with
 holomorphic functions with compact support. For this let consider
 the function ${\Ker}^\nu_{\Gamma,\chi}(z,w)$ defined on
 $\C^n\times\C^n$ by the convergent series in $\mathcal{C}^\infty(\C^n\times\C^n)$
 \begin{equation}\label{RepKer}
 {\Ker}^\nu_{\Gamma,\chi}(z,w) := \Big(\frac\nu \pi\Big)^{n}e^{\nu\scal{z,w}} \sum_{\gamma\in \Gamma}
\chi(\gamma)e^{-\frac \nu
2|\gamma|^2+\nu(\scal{z,\gamma}-\overline{\scal{w,\gamma}})}.
\end{equation}
Then, we state

 \begin{theorem}\label{ThmRepKer} Suppose that the condition $(RDQ)$
 is satisfied and let ${\Ker}^\nu_{\Gamma,\chi}(z,w)$ be the function defined by
 (\ref{RepKer}). Then

    \begin{itemize}
          \item[i)] For every $z,w\in \C^n$, we have ${\Ker}^\nu_{\Gamma,\chi}(z,w)
 =\overline{{\Ker}^\nu_{\Gamma,\chi}(w,z)} $.
          \item[ii)]For every $\gamma_1,\gamma_2\in \Gamma$ and $z,w\in \C^n$, we
have \begin{equation}
{\Ker}^\nu_{\Gamma,\chi}(z+\gamma_1,w+\gamma_2) =
\chi(\gamma_1)e^{\frac \nu 2|\gamma_1|^2+\nu\scal{z,\gamma_1}}
{\Ker}^\nu_{\Gamma,\chi}(z,w)\overline{\chi(\gamma_2)}e^{\frac \nu
2|\gamma_2|^2+\nu\overline{\scal{w,\gamma_2}}}.
\label{GammaInvRepKer}
\end{equation}
In particular the function $z
\longmapsto{\Ker}^\nu_{\Gamma,\chi}(z,w)$ belongs to
${\mathcal{O}}^\nu_{\Gamma,\chi}(\C^n)$ for every fixed $w\in\C^n$.
          \item[iii)]  For $\mbox{vol}(\C^n/\Gamma)$ denoting the Lebesgue volume of a
fundamental domain of the lattice $\Gamma$, we have
 \[\int_{\Lambda(\Gamma)}
{\Ker}^\nu_{\Gamma,\chi}(z,z)e^{-\nu|z|^2}dm(z) = ({\nu}/{\pi})^n
\mbox{vol}(\C^n/\Gamma).\]
          \item[iv)] For every $g\in {\mathcal{O}}^\nu_{\Gamma,\chi}(\C^n)$, we have
\[ g(z) = \int_{\Lambda(\Gamma)} {\Ker}^\nu_{\Gamma,\chi}(z,w)g(w)
e^{-\nu|w|^2} dm(w).\]
    \end{itemize}
 \end{theorem}

 \noindent{\it Proof.} i) is easy to check. Indeed we use the fact
 that $\overline{\chi(\gamma)}=\chi(-\gamma)$, for the $(RDQ)$ condition being satisfied, and next the
 change $-\gamma$ by $\gamma$ in the involved summation (\ref{RepKer}).

For the assertion ii),  one gets from (\ref{RepKer}) that
\[ {\Ker}^\nu_{\Gamma,\chi}(z+\gamma,w) =({\nu}/{\pi})^n e^{\nu\scal{z,w}}
\sum_{\gamma'\in \Gamma} \chi(\gamma') e^{\nu\scal{\gamma,\gamma'}}
e^{-\frac \nu
2|\gamma'|^2+\nu(\scal{z,\gamma'}-\overline{\scal{w,\gamma'-\gamma}})}.\]
Next, by making the change $\gamma"=\gamma'-\gamma$ and using the
$(RDQ)$ condition, it follows
\[ {\Ker}^\nu_{\Gamma,\chi}(z+\gamma,w) =
\chi(\gamma)e^{\frac \nu 2|\gamma|^2+\nu\scal{z,\gamma}}
{\Ker}^\nu_{\Gamma,\chi}(z,w); \qquad \gamma\in \Gamma.\]
 Hence we obtain $(\ref{GammaInvRepKer})$ thanks to i).
 Therefore, the function ${\Ker}^\nu_{\Gamma,\chi}(z,w)$
belongs to ${\mathcal{O}}^\nu_{\Gamma,\chi}(\C^n)$ for being
holomorphic function.

For iii), we make use of (\ref{RepKer}) again, to get
\begin{align*}  \int_{\Lambda(\Gamma)}
{\Ker}^{\nu}_{\Gamma,\chi}(z,z)e^{-\nu|z|^2}dm(z)
 = ({\nu}/{\pi})^n\sum_{\gamma\in \Gamma} \chi(\gamma) e^{-\frac \nu
2|\gamma|^2} \bigg(\int_{\Lambda(\Gamma)}e^{2i\nu\omega(z,\gamma)}
dm(z)\bigg ).
\end{align*}
Now, since the condition $(RDQ)$ is satisfied, we see that the
function $z\mapsto e^{2i\nu\omega(z,\gamma)} $ is $\Gamma$-periodic
for every fixed $\gamma\in \Gamma$ and therefore $S_\gamma$;
$\gamma\in \Gamma$, where $S_\gamma(z):=e^{2i\nu \omega(z,\gamma)}$,
define a group character on $\C^n/\Gamma$. Hence, we have (see
\cite[page 3480]{Bump}, but a direct proof is presented hereafter):
\begin{lemma}\label{Integral0}
Assume that the $(RDQ)$ condition is verified. Then, for every
 $\gamma\in \Gamma\setminus\set{0}$, we have
\begin{equation}
\int_{\Lambda(\Gamma)} e^{2i\nu\omega(w,\gamma)} dm(w) = 0.
\label{AverageCell}\end{equation}
\end{lemma}
\noindent Thus, it follows
\begin{align*}  \int_{\Lambda(\Gamma)}
{\Ker}^{\nu}_{\Gamma,\chi}(z,z)e^{-\nu|z|^2}dm(z) =
({\nu}/{\pi})^n\chi(0) \mbox{vol}(\Lambda(\Gamma))=
({\nu}/{\pi})^n\mbox{vol}(\Lambda(\Gamma)).
\end{align*}

 The proof of iv)  relies essentially on the following
\begin{lemma}\label{KeyLem}  For any given holomorphic function $g\in {\mathcal{O}(\C^n)}$
satisfying the growth condition
 $|g(z)| \leq C e^{\frac \nu 2 |z|^2}$,
we have the following reproducing formula
\begin{equation} g(z) =  ({\nu}/{\pi})^n \int_{\C^n}  e^{\nu\scal{z,w}} g(w)
e^{-\nu|w|^2} dm(w). \label{KeyLemma}\end{equation}
 \end{lemma}
\noindent Then, since for every $g\in
{\mathcal{O}}^\nu_{\Gamma,\chi}(\C^n)$, there exists certain
constant $C\geq 0$ such that $|g(z)| \leq C e^{\frac \nu 2 |z|^2}$,
one can apply (\ref{KeyLemma}) to have
\[
g(z) =   ({\nu}/{\pi})^n \int_{\C^n} e^{\nu\scal{z,w}} g(w)
e^{-\nu|w|^2} dm(w).
\]
Next, by writing $\C^n$ as disjoint union of $ \gamma+
\Lambda(\Gamma)$, for varying $\gamma\in\Gamma$ and using the fact
that the function $g$ satisfies the functional equation
\[ g(w+\gamma)=\chi(\gamma)e^{\frac\nu 2|\gamma|^2+\nu\scal{w,\gamma}}g(w) ,\] we get
\begin{align*}
g(z)  &= ({\nu}/{\pi})^n\sum_{\gamma\in
\Gamma}\int_{\Lambda(\Gamma)} e^{\nu\scal{z,w+\gamma}} g(w+\gamma)
e^{-\nu|w+\gamma|^2} dm(w)
 \\
 &= \int_{\Lambda(\Gamma)} \Big[({\nu}/{\pi})^n\sum_{\gamma\in \Gamma}
\chi(\gamma)e^{\nu\scal{z,w+\gamma}}e^{-\frac \nu
2|\gamma|^2-\nu\overline{\scal{w,\gamma}}}\Big] g(w) e^{-\nu|w|^2}
dm(w).
\end{align*}
Finally by definition of ${\Ker}^\nu_{\Gamma,\chi}(z,w)$ given by
(\ref{RepKer}), we conclude that
 \[ g(z)=
\int_{\Lambda(\Gamma)} {\Ker}^\nu_{\Gamma,\chi}(z,w) g(w)
e^{-\nu|w|^2} dm(w).   \eqno{\mbox{\fin}}
\]

 \noindent{\it Proof of Lemma \ref{Integral0}.} Fix
$\gamma\in\Gamma$ such that $\gamma\ne 0$ and note that $z \to
\omega(z,\gamma)$ is a nonzero function on a given fundamental
domain $\Lambda(\Gamma)$. Let $u_1,u_2, \cdots , u_{2n} \in \Gamma$
be a basis of $\Lambda(\Gamma)$. For every fixed $z\in
\Lambda(\Gamma)$, we write $z=t_1u_1+t_2u_2+ \cdots +t_{2n}u_{2n}$
with $t_j \in[0,1]$. By the $G$-invariance of $dm$, we get
\begin{align*} \int_{\Lambda(\Gamma)}e^{2i\nu \omega(z,\gamma)}
dm(z) &=  vol(\Lambda(\Gamma))\prod_{j=1}^{2n} \int_0^1 e^{2i\nu
t_j\omega(u_j,\gamma)} dt_j
\\
&=vol(\Lambda(\Gamma))\prod_{{\tiny{\begin{array}{c} j=1 \\
\omega(u_j,\gamma)\ne 0\end{array}}} }^{2n} \int_0^1 e^{2i\nu
t_j\omega(u_j,\gamma)} dt_j
\\&= vol(\Lambda(\Gamma))\prod_{{\tiny{\begin{array}{c} j=1 \\ \omega(u_j,\gamma)\ne
0\end{array}}} }^{2n} \frac{\Big( e^{2i\nu
\omega(u_j,\gamma)}-1\Big)}{2i\nu \omega(u_j,\gamma)}.
\end{align*}
But, since  $ e^{2i\nu \omega(u_j,\gamma)} =1 $ for the triplet
$(\nu,\Gamma,\chi)$ satisfying the $(RDQ)$ condition, we get
\[\int_{\Lambda(\Gamma)}e^{2i\nu \omega(z,\gamma)} dm(z) =0.\eqno{\mbox{\fin}}
\]

\noindent{\it Proof of Lemma \ref{KeyLem}.} Let $g\in
{\mathcal{O}(\C^n)}$ such that
 $|g(z)| \leq C e^{\frac \nu 2 |z|^2}$ and consider the function $ g_\varepsilon
(z):= g(\varepsilon z)$; $z\in \C^n$, for every given $\varepsilon$;
$0<\varepsilon <1$. Then, clearly $g_\varepsilon$ is holomorphic and
satisfies $|g_\varepsilon(z)| \leq C e^{\frac \nu 2 \varepsilon
^2|z|^2}.$ Furthermore, we have
\[ \int_{\C^n} |g_\varepsilon(z)|^2 e^{-\nu|z|^2} dm(z) \leq C^2
\int_{\C^n} e^{-\nu(1-\varepsilon^2)|z|^2} dm(z)< +\infty\] and
therefore $g_\varepsilon$ belongs to the Bargmann-Fock space $
{\mathcal{B}}^{2,\nu}(\C^n)$ on $\C^n$,
\begin{equation} {\mathcal{B}}^{2,\nu}(\C^n)=
                     \set{g  \in {\mathcal{O}(\C^n)} ;
                  \quad \int_{\C^n} |g(z)|^2 e^{-\nu
                   |z|^2}dm(z) < +\infty}.\label{BargmannSpace}
 \end{equation}
Hence, for every $\varepsilon$ with $0<\varepsilon <1$, we have
\[ g_\varepsilon(z)
=  ({\nu}/{\pi})^n\int_{\C^n}e^{\nu\scal{z,w}} g_\varepsilon(w)
e^{-\nu|w|^2} dm(w).\]
 Finally by tending $\varepsilon$ to $1$ and applying
the dominated convergence theorem, we get
\[ g(z)
=  ({\nu}/{\pi})^n\int_{\C^n}e^{\nu\scal{z,w}}  g(w) e^{-\nu|w|^2}
dm(w). \eqno{\mbox{\fin}}\]

As consequence of the previous theorem, we state the following
result for the space $ {\mathcal{O}}^\nu_{\Gamma,\chi}(\C^n)$
endowed with the norm $||\cdot||_\Gamma$ associated with
(\ref{HSPGamma2}).

\begin{corollary} \label{LastRepKer}
Let the condition $(RDQ)$ to be satisfied. Then

    \begin{itemize}
          \item[i)]  ${\mathcal{O}}^\nu_{\Gamma,\chi}(\C^n)$ is a nonzero space.
          \item[ii)]  ${\mathcal{O}}^\nu_{\Gamma,\chi}(\C^n)$ is a reproducing Hilbert
space whose the reproducing kernel is given by (\ref{RepKer}), i.e.,
\begin{equation}\label{RepKer2}
 {\Ker}^\nu_{\Gamma,\chi}(z,w) := \Big(\frac\nu \pi\Big)^{n}e^{\nu\scal{z,w}} \sum_{\gamma\in \Gamma}
\chi(\gamma)e^{-\frac \nu
2|\gamma|^2+\nu(\scal{z,\gamma}-\overline{\scal{w,\gamma}})}.
\end{equation}
          \item[iii)]  ${\mathcal{O}}^\nu_{\Gamma,\chi}(\C^n)$ is a finite dimensional
space whose dimension is \[\dim
{\mathcal{O}}^\nu_{\Gamma,\chi}(\C^n) =({\nu}/{\pi})^n
\mbox{vol}(\Lambda(\Gamma)).\]
    \end{itemize}
\end{corollary}

\noindent{\it Proof.} For i), we see from iii) of Theorem
\ref{ThmRepKer}  that ${\Ker}^\nu_{\Gamma,\chi}(z,w)$ is a non
vanishing function on $\C^n\times \C^n$. Hence
${\Ker}^\nu_{\Gamma,\chi}(z_0,w_0)\ne 0$ for some $(z_0,w_0)\in
\C^n\times \C^n$. Therefore the function $z \mapsto
{\Ker}^\nu_{\Gamma,\chi}(z,w_0)$ is a nonzero function that belongs
to ${\mathcal{O}}^\nu_{\Gamma,\chi}(\C^n)$ by ii) of Theorem
\ref{ThmRepKer}.

To prove ii), we apply the Cauchy-Schwartz inequality to iv) of
Theorem \ref{ThmRepKer}. Thus, we see that for every $g\in
{\mathcal{O}}^\nu_{\Gamma,\chi}(\C^n)$, we have
\begin{equation}|g(z)| \leq   \Big(\int_{\Lambda(\Gamma)}
|{\Ker}^\nu_{\Gamma,\chi}(z,w)|^2 e^{-\nu|w|^2} dm(w)\Big)^{1/2}
||g||_\Gamma.\label{Ineq}\end{equation} Then for any given bounded
domain $\Omega\subset\C^n$, we get
\begin{equation}|g(z)| \leq   C_\Omega
||g||_\Gamma, \quad z\in
\overline{\Omega},\label{Ineq2}\end{equation} for certain constant
$C_\Omega$.
 Therefore, any Cauchy sequence $g_j$ in
${\mathcal{O}}^\nu_{\Gamma,\chi}(\C^n)$ for the norm
$||\cdot||_\Gamma$, is also a uniformly Cauchy sequence on any
compact set of $\C^n$ and then converges to a holomorphic function
$g$ on $\C^n$. Next, since all $g_j$ satisfies the functional
equation
\[g_j(z+\gamma) = \chi(\gamma) e^{\frac \nu 2
|\gamma|^2+\nu\scal{z,\gamma}}g_j(z),\] it follows that $g$ belongs
to ${\mathcal{O}}^\nu_{\Gamma,\chi}(\C^n)$.

 This shows that the nonzero space ${\mathcal{O}}^\nu_{\Gamma,\chi}(\C^n)$ is in fact a Hilbert space
  for the norm $||\cdot ||_\Gamma$. Furthermore, the use of (\ref{Ineq2})
 infers that the
 evaluation map $\delta_{z}:
{\mathcal{O}}^\nu_{\Gamma,\chi}(\C^n) \longrightarrow \C$ given by $
\delta_{z} g := g(z) ,$ is continuous for every fixed $z\in \C^n$.
Hence ${\mathcal{O}}^\nu_{\Gamma,\chi}(\C^n)$ possesses a unique
reproducing kernel function $\tilde K^\nu_{\Gamma,\chi}(z,w)$. But
in view of Theorem 3.7 [mainly i), ii) and iv)], we deduce that
 $\tilde K^\nu_{\Gamma,\chi}(z,w)$ is exactly the function
${\Ker}^\nu_{\Gamma,\chi}(z,w)$ given through (\ref{RepKer2}).

The dimension of the Hilbert
  space ${\mathcal{O}}^\nu_{\Gamma,\chi}(\C^n)$ can be calculated by integrating its
  reproducing kernel ${\Ker}^\nu_{\Gamma,\chi}(z,w)$ along the diagonal, that is
\[ \dim {\mathcal{O}}^\nu_{\Gamma,\chi}(\C^n) =
\int_{\Lambda(\Gamma)}{\Ker}^\nu_{\Gamma,\chi}(z,z)
 e^{-\nu |z|^2} dm(z).
\]
Whence, in light of ii) of Theorem \ref{ThmRepKer}, it follows \[\dim
{\mathcal{O}}^\nu_{\Gamma,\chi}(\C^n)= ({\nu}/{\pi})^n
 \mbox{vol}(\Lambda(\Gamma))\]
 and hence the proof of the theorem is completed.
 \fin

\begin{remark}\label{RemPeriod}   The space ${\mathcal{O}}^\nu_{\Gamma,\chi}(\C^n)$ is of interest in
itself for it linked in somehow to the classical Bargmann-Fock space $
{\mathcal{B}}^{2,\nu}(\C^n)$.
 Indeed, according to the explicit expression (\ref{RepKer}), the
 reproducing kernel ${\Ker}^\nu_{\Gamma,\chi}$
 of the space ${\mathcal{O}}^\nu_{\Gamma,\chi}(\C^n)$ appears then as
 $(\Gamma,\chi)$-periodization, with respect to the automorphy
 factor $\chi(\gamma)e^{\frac \nu
2|\gamma|^2+\nu{\scal{w,\gamma}}}$,
 of the reproducing kernel $({\nu}/{\pi})^ne^{\nu{\scal{w,\gamma}}}$ of the Bargmann-Fock
space $\mathcal{B}^{2,\nu}(\C^n)$ (\ref{BargmannSpace}).
\end{remark}

In the following table, we summarize some basic properties related
to the spaces ${\mathcal{F}}^\nu_{\Gamma,\chi}$,
${\mathcal{G}}^\nu_{\Gamma,\chi}$ and
${\mathcal{O}}^\nu_{\Gamma,\chi}(\C^n)$.

\vspace*{.3cm}

\begin{center}
 {\small{\begin{tabular}{|c|c|c|c|c|}
  \hline
  & Related items to ${\mathcal{F}}^\nu_{\Gamma,\chi}$ &  $\mathfrak{G}$   &  Related items to
  ${\mathcal{G}}^\nu_{\Gamma,\chi}$ \\
  \hline \hline
  Automorphy factor     &
    $\chi(\gamma) e^{i\nu\omega(z,\gamma)}$  &     & $\chi(\gamma) e^{\frac \nu 2|\gamma|^2+\nu\scal{z,\gamma}}$ \\ \hline
  Functional equation    &
 $f(z+\gamma) = \chi(\gamma) e^{i\nu\omega(z,\gamma)}f(z)$   &      & $g(z+\gamma) = \chi(\gamma) e^{\frac \nu 2
                                                                 |\gamma|^2+\nu\scal{z,\gamma}}g(z)$  \\\hline
  Growth condition       &
   $ |f(z)|\leq C$ &      & $|g(z)|\leq Ce^{\frac \nu 2 |z|^2}$ \\\hline
     Natural subspaces      &
  $ \mathfrak{G}^{-1}[ {\mathcal{O}}^\nu_{\Gamma,\chi}(\C^n)] $ \quad \fbox{?} &   $\cong$   & $ {\mathcal{O}}^\nu_{\Gamma,\chi}(\C^n)$ \\\hline
     Dimension formulas      &
  $ \dim\mathfrak{G}^{-1}[ {\mathcal{O}}^\nu_{\Gamma,\chi}(\C^n)] $ & =
   & $ \dim{\mathcal{O}}^\nu_{\Gamma,\chi}(\C^n)$\\\hline
  Scalar product         &
  $\scal{f_1, f_2}_\Gamma$ &     &
  $\scal{\scal{g_1, g_2}}_\Gamma$ \\\hline
  Hilbert structure       &
    $L^{2,\nu}_{\Gamma,\chi}=\overline{{\mathcal{F}}^\nu_{\Gamma,\chi}}^{\scal{, }_\Gamma}$ &
    $\cong$  & $
    B^{2,\nu}_{\Gamma,\chi}=\overline{{\mathcal{G}}^\nu_{\Gamma,\chi}}
     ^{\scal{\scal{\cdot , \cdot }}_\Gamma}$ \\\hline
  Differential operator  &
   $\La^\nu_{\Gamma,\chi}$ &     & $\Delta^\nu_{\Gamma,\chi} $
   \\\hline
     Particular Eigenspaces      &
  $  \ker (\La^\nu_{\Gamma,\chi}- n\nu) ={\mathcal{E}}^\nu_{\Gamma,\chi}(0) $ &  $\cong$   & $ \ker \Delta^\nu_{\Gamma,\chi} $
   \\\hline
     Dimension Formulas      &
  $  \dim{\mathcal{E}}^\nu_{\Gamma,\chi}(0)  $ \quad \fbox{?} &   $=$
    & $ \dim\ker \Delta^\nu_{\Gamma,\chi} $\quad \fbox{?}\\ \hline
  \hline
\end{tabular}}}
\end{center}

\begin{remark}
The answers to the three question marks in the above table are
included in the main theorem. and iii) of Corollary
\ref{LastRepKer}. In fact, the space
 $ \mathfrak{G}^{-1}[{\mathcal{O}}^\nu_{\Gamma,\chi}(\C^n)]$
 represents in the ${\mathcal{F}}^\nu_{\Gamma,\chi}$ picture the
 eigenspace of $\La^\nu_{\Gamma,\chi}$ associated with its lowest
 eigenvalue. Thus $$ {\mathcal{O}}^\nu_{\Gamma,\chi}(\C^n)\cong
 \mathfrak{G}^{-1}[{\mathcal{O}}^\nu_{\Gamma,\chi}(\C^n)]
 \equiv {\mathcal{E}}^\nu_{\Gamma,\chi}(0) \cong \ker
 \Delta^\nu_{\Gamma,\chi}.$$
\end{remark}

 \section{General properties of the
      $(\Gamma,\chi)$-automorphic kernel functions of magnitude $\nu>0$}

Here we reconsider the action of the semi-direct group  $G=U(n)
\rtimes \C^n$ on $L^2(\C^n;dm)$ given through the unitary
transformations (\ref{MT01}), $ [T^\nu_{g}f](z):=j_\nu(g,z) f(g.z),
$ where $j_\nu(g,z)$ $= e^{i\nu\omega(z,g^{-1}.0)}$.

\begin{definition} A given $\mathcal{C}^\infty$ function $\mathcal{K}(z,w)$ on $\C^n\times \C^n$ is said to be
$G$-invariant if it verifies the following property
\begin{equation} {\mathcal{K}}(g.z,g.w)  =
 \overline{j_{\nu}(g,z)}
{\mathcal{K}} (z,w) j_{\nu}(g,w) .\label{IPKerFun}\end{equation}
\end{definition}
 Thus, one can see that a given kernel function $\mathcal{K}(z,w)$ on $\C^n\times \C^n$ is
$G$-invariant if and only if it is of the form \begin{equation}
{\mathcal{K}}(z,w) = e^{i{\nu}\omega(z,w)}
 Q^\nu(|z-w|)\label{CanoFormInvKer}
\end{equation}
for certain function $Q^\nu$ defined on the positive real line.
Throughout this section, we will suppose that the
 involved $Q^\nu$ is rapidly decreasing function on $[0,+\infty)$
 and that the $(RDQ)$ condition is satisfied.
 Thus, we define
${\mathcal{K}}^{\nu}_{\Gamma,\chi}$ to be
 the function on $\C^n\times\C^n$ given by the following convergent
 series
 \begin{equation}
\label{IPKerFunForm} {\mathcal{K}}^{\nu}_{\Gamma,\chi} (z,w) =
e^{i\nu\omega(z,w)}
 \sum_{\gamma\in \Gamma}
\chi(\gamma)e^{i\nu\omega(z+w,\gamma)} Q^\nu(|z-w-\gamma|).
\end{equation}
 The function
$ {\mathcal{K}}^{\nu}_{\Gamma,\chi}(z,w)$ is in fact the
($\Gamma,\chi$)-periodization (\`a la Poincar\'e) of the appropriate
function ${\mathcal{K}}_z(w):= {\mathcal{K}}(z,w)$ and can be
rewritten in the following variant forms
\begin{align}
 {\mathcal{K}}^{\nu}_{\Gamma,\chi}(z,w)
&= \sum_{\gamma\in \Gamma} \chi(\gamma) T^{-\nu}_{\gamma}
[{\mathcal{K}}_z](w)\label{GammaRepKer21}\\
&=  \sum_{\gamma\in \Gamma} \chi(\gamma) T^{\nu}_{-\gamma}
[{\mathcal{K}}(\xi,w)]|_{\xi=z} \label{GammaRepKer12}
\end{align}
 Therefore, 
 using the $(RDQ)$ condition, it can be shown that
${\mathcal{K}}^{\nu}_{\Gamma,\chi} (z,w)$ satisfies the following
$\Gamma$-bi-invariant property
 \begin{equation} \label{IPKerFunGammaBi}
 {\mathcal{K}}^{\nu}_{\Gamma,\chi}
(z+\gamma,w+\gamma') = \chi(\gamma)
 \overline{j_{\nu}(\gamma,z)} {\mathcal{K}}^{\nu}_{\Gamma,\chi}
 (z,w) \overline{\chi(\gamma')}j_{\nu}(\gamma',w)
\end{equation}
for every $\gamma,\gamma'\in \Gamma$, and in particular it is
$\Gamma$-invariant. Moreover, we have
$\overline{{\mathcal{K}}^{\nu}_{\Gamma,\chi} (z,w)} =
{\mathcal{K}}^{\nu}_{\Gamma,\chi} (w,z)$ if and only if the function
$Q^\nu$ is assumed to be, in addition, a real valued function.

\begin{definition} We call ${\mathcal{K}}^{\nu}_{\Gamma,\chi} (z,w)$, given
by (\ref{IPKerFunForm}),   the $(\Gamma,\chi)$-automorphic kernel
function  on $\C^n\times \C^n$ of magnitude $\nu>0$ associated with
the $G$-invariant kernel function ${\mathcal{K}}(z,w)$.
\end{definition}

 Now, let denote by $\K$ the integral operator acting
 on the Hilbert space $L^2(\C^n;dm)$
  by
 \begin{align}\Big[\K(\varphi)\Big](z) &=
\int_{\C^n} {\mathcal{K}}(z,w) \varphi(w) dm(w)
\label{IntOpKerFun}\\
&=  \int_{\C^n}e^{i{\nu}\omega(z,w)}
 Q^\nu(|z-w|)\varphi(w) dm(w)\label{IntOpKerFun2}\end{align}
which is well defined for $ Q^\nu$ being assumed to be a rapidly
decreasing function. Also, let denote by $\K^{\nu}_{\Gamma,\chi}$
the integral operator  associated with the $(\Gamma,\chi)$-automorphic
kernel function ${\mathcal{K}}^{\nu}_{\Gamma,\chi}$ and acting on
the Hilbert space $L^{2,\nu}_{\Gamma,\chi}$ by
 \begin{equation}\Big[\K^{\nu}_{\Gamma,\chi}(\psi)\Big](z) =
\int_{\Lambda(\Gamma)} {\mathcal{K}}^{\nu}_{\Gamma,\chi}(z,w)
\psi(w) dm(w),\label{IntOpKerFunGamma}.\end{equation}

At once, since $Q^\nu$ is rapidly decreasing, we show that the
integral operator $\K^{\nu}_{\Gamma,\chi}$ is of trace. More
precisely, we have the following result (whose the proof is exactly
the same as the one provided for ii) of Theorem \ref{ThmRepKer}).
\begin{proposition}\label{PropTrace} Under the $(RDQ)$ condition and the
assumption that $Q^\nu$ is rapidly decreasing,
 the trace of the integral operator
$ \K^{\nu}_{\Gamma,\chi}$ can be given by \begin{equation}
\label{Trace} Trace(\K^{\nu}_{\Gamma,\chi})=\int_{\Lambda(\Gamma)}
{\mathcal{K}}^{\nu}_{\Gamma,\chi}(z,z) dm(z) = Q^\nu(0)
\rm{vol}(\Lambda(\Gamma)).\end{equation}
\end{proposition}
  As immediate consequence, we have
\begin{corollary} \label{W0}
If the function $Q^\nu$ verifies $Q^\nu(0)\ne 0$, then there exists
$w_0 \in \C^n$ such that the function $ z \longmapsto
{\mathcal{K}}^{\nu}_{\Gamma,\chi}(z,w_0)$ is  nonzero function on
$\C^n$.
\end{corollary}

The relationship between $\K$ and $\K^{\nu}_{\Gamma,\chi}$ is given
by the following
\begin{lemma} \label{RepKerResLem}
For every $\psi \in \mathcal{F}^{\nu}_{\Gamma,\chi}$, we have
\[ [\K(\psi)](z) := \int_{\C^n} {\mathcal{K}}(z,w) \psi(w) dm(w)= 
\int_{\Lambda(\Gamma)}
\mathcal{K}^{\nu}_{\Gamma,\chi}(z,w)\psi(w)dm(w)\] which means that
$\K_{|_{\mathcal{F}^{\nu}_{\Gamma,\chi}}} = \K^{\nu}_{\Gamma,\chi}.
$
\end{lemma}

\noindent{\it Proof.} Writing
 $\C^n$ as disjoint union of $ \gamma+\Lambda(\Gamma)$, $\gamma\in\Gamma$, and using
   the fact that the Lebesgue measure
 $dm$ is $G$-invariant as well as that any arbitrary function $\psi$ in
$\mathcal{F}^{\nu}_{\Gamma,\chi}$ satisfies
$\psi(w+\gamma )=\chi(\gamma) \overline{j_\nu(\gamma,w)}
\psi(w) $ for every $\gamma\in\Gamma$ and every $w\in \C^n$. \fin

\vspace*{.2cm} The above property allows $\K^{\nu}_{\Gamma,\chi}$ to
inherit some useful properties of $\K$. For instance, we have the
following commutations:

\begin{proposition} \label{KerFunCommu} \quad

    \begin{itemize}
          \item[i)]  For every $g\in G$, we have $
 T^\nu_{g}\K = \K T^\nu_{g}. $
          \item[ii)]  For every $\gamma\in \Gamma$, we have
$T^\nu_{\gamma}\K^{\nu}_{\Gamma,\chi}= \K^{\nu}_{\Gamma,\chi}
T^\nu_{\gamma}.$
          \item[iii)]  The  Landau Hamiltonian $\La^\nu$
 commutes with both integral operators $\K$ and $\K^{\nu}_{\Gamma,\chi}$, i.e.,
\begin{equation} \label{CommRule3}
\La^\nu\K = \K  \La^\nu \quad \mbox{ and } \quad
\La^\nu\K^{\nu}_{\Gamma,\chi}= \K^{\nu}_{\Gamma,\chi}
\La^\nu.\end{equation}
    \end{itemize}
\end{proposition}

\noindent{\it Proof.} i) is easy to check. Indeed, using the
invariance property (\ref{IPKerFun}) that we can rewrite also as
\begin{equation} j_{\nu}(g,z){\mathcal{K}}(g.z,w) =
 {\mathcal{K}}(z,g^{-1}.w)
 \overline{j_{\nu}(g^{-1},w)}
,\label{IPKerFun2}\end{equation}
 it follows that
 \[
 T^\nu_{g}\big[\K(\varphi) \big](z) \stackrel{(\ref{IPKerFun2})}{=} \int_{\C^n}
                                    {\mathcal{K}}(z,g^{-1}.w) \overline{j_\nu(g^{-1},w)} \varphi (w) dm(w).
\]
Next, by making use of the change $w'=g^{-1}.w$, we conclude that
 \[
 T^\nu_{g}\big[\K(\varphi) \big](z) =  \int_{\C^n}  {\mathcal{K}}(g,w')  [T^\nu_{g}(\varphi)] (w') dm(w')
                                    =  \K \big[T^\nu_{g}(\varphi)\big] (z)
\]

 ii) is an immediate consequence of Lemma \ref{RepKerResLem}
combined with i) above, keeping in mind the fact that $
T^\nu_{\gamma}\psi$ belongs to $\mathcal{F}^{\nu}_{\Gamma,\chi}$ if
$\psi\in\mathcal{F}^{\nu}_{\Gamma,\chi}$.

 For iii), we begin by noting that if we use the notation $\La^\nu_u$
 to mean that the derivation is taken w.r.t. the complex variable $u$,
 then we have
\begin{equation}
\label{LZW} \La^{\nu}_z{\mathcal{K}}(z,w)=
\La^{-\nu}_w{\mathcal{K}}(z,w).
\end{equation}
 Put ${\mathcal{K}_0}(\xi):={\mathcal{K}}(\xi,0)$. Hence in view of the
fact that for every  $g_z,g_w\in G$ such that $g_z.z=0$ and
$g_w.w=0$, we have
\begin{equation}\label{KerTran}
{\mathcal{K}}(z,w) =
\big[T^\nu_{g_w}\big({\mathcal{K}_0}(\xi)\big)\big]|_{\xi=z} =
\big[T^{-\nu}_{g_z}\big({\mathcal{K}_0}(\xi)\big)\big]|_{\xi=w},\end{equation}
 together with the fact that $\La^\nu_z$ and
$T^\nu_{g_w}$ commute, we get
\begin{equation*}[\La^\nu_z{\mathcal{K}}(z,w)] =
\La^\nu_\xi[T^\nu_{g_w}{\mathcal{K}_0}(\xi)]|_{\xi=z} =
T^\nu_{g_w}[\La^\nu_\xi{\mathcal{K}_0}(\xi)]|_{\xi=z} =
T^\nu_{g_w}[\La^{-\nu}_\xi{\mathcal{K}_0}(\xi)]|_{\xi=z}.
\end{equation*}
The last equality follows from the facts that $\xi\to
{\mathcal{K}_0}(\xi)$ is radial and the operators $\La^{-\nu}$ and
$\La^\nu$ have the same radial parts. Next, using the observation
\[ [T^\nu_{g_w}\varphi](z) =
[T^{-\nu}_{g_z}\varphi](w),\] for any radial function $\varphi$, we
conclude
\begin{align*}
[\La^\nu_z{\mathcal{K}}(z,w)]&=
T^{-\nu}_{g_z}[\La^{-\nu}_\xi{\mathcal{K}_0}(\xi)]|_{\xi=w}=
\La^{-\nu}_\xi[T^{-\nu}_{g_z}{\mathcal{K}_0}(\xi)]|_{\xi=w}
\stackrel{(\ref{KerTran})}{=} [\La^{-\nu}_w{\mathcal{K}}(z,w)].
\end{align*}
Therefore, it follows
\[
\La^\nu[\K f](z)= \int_{\C^n}[\La^{-\nu}_w{\mathcal{K}}(z,w)]f(w)
 dm(w)
 \]
 for every $\mathcal{C}^\infty$ compactly supported function $f$.
 Finally, integration by parts yields
 \[
\La^\nu [\K f](z) = \K [\La^\nu f](z).
 \]
The commutation $\La^\nu\K^{\nu}_{\Gamma,\chi}=
\K^{\nu}_{\Gamma,\chi} \La^\nu$ follows easily from the previous one
using Lemma \ref{RepKerResLem} together with observation that
$\La^\nu f \in \mathcal{F}^{\nu}_{\Gamma,\chi}$ for $f\in
\mathcal{F}^{\nu}_{\Gamma,\chi}$. \fin
\begin{remark} It can be shown that integral operators $\K$
satisfying the commutation rule $T^\nu_{g}\K = \K T^\nu_{g}$ are
those whose kernel function $\mathcal{K}(z,w)$ is $G$-invariant and
so is of the form
\begin{equation}\label{RedForm} \mathcal{K}(z,w)= e^{i\nu\omega(z,w)} Q^\nu(|z-w|).\end{equation}
\end{remark}

 To state the next result, let $g\in G$ and denote by
$\mathbf{{{\Av}_g}}$ the averaging operator acting on $L^2(\C^n;dm)$
by
\begin{equation}[ {{\Av}_g}(\psi)](z):=
  \int_{U(n)}\Big[T^\nu_{g k} (\psi)\Big](z) dk, \label{average0}
\end{equation} where $dk$ is the normalized
 Haar measure of ~$U(n)$. Then, we state
\begin{proposition}\label{PropSelberg1}
The integral operator  $\K$ associated with the $G$-invariant kernel
function ${\mathcal{K}}$ commutes with the averaging operator
${{\Av}_g}$ as defined in (\ref{average0}). That is for every $g\in
G$ and every $\psi \in L^2(\C^n;dm)$, we have
\begin{equation}
 {{\Av}_g} [\K (\psi)](z)
  = \K[{{\Av}_g}(\psi)](z). \label{IntOpAveraging}
\end{equation}
\end{proposition}

\noindent{\it Proof.} By definition, we have
\[
{{\Av}_g} [\K (\psi)](z) =\int_{U(n)}T^\nu_{g k}[\K (\psi)](z) dk.
\]
Next, by applying i) of Proposition \ref{KerFunCommu},  it follows
\begin{align*}
{{\Av}_g} [\K (\psi)](z) &= \int_{U(n)}\K[T^\nu_{g
k}(\psi)](z) dk\\
&=   \int_{U(n)}\int_{\C^n}{\mathcal{K}}(z,w)[T^\nu_{g k}(\psi)](w)
 dm(w) dk.
 \\
&=   \int_{\C^n}{\mathcal{K}}(z,w)\Big( \int_{U(n)} [T^\nu_{g
k}(\psi)](w) dk \Big)
 dm(w) \\
 &=   \int_{\C^n}{\mathcal{K}}(z,w)[{{\Av}_g}(\psi)](w)dm(w)
\\
&=
 \K [{{\Av}_g}(\psi)](z). 
\end{align*}
\fin

\noindent Therefore, one can establish the following properties of
the above averaging operator ${{\Av}_g}$. Namely, we have

\begin{proposition} \label{AveragingFunction}
Let $\psi\in L^2(\C^n;dm)$ and ${{\Av}_g}(\psi)$ its averaging
  as defined by (\ref{average0}) for given $g\in G$. Then,

    \begin{itemize}
          \item[i)]  ${{\Av}_g}\psi$ is a radial function for every $g\in G$.
          \item[ii)] If $\psi$ is  bounded, then ${{\Av}_g}\psi$ is bounded for every
                   $g\in G$.
          \item[iii)] The averaging operator ${{\Av}_g}$ commutes with the Hamiltonian $\La^\nu$.
          \item[iv)]  If $\psi$ is a bounded nonzero solution of the differential
equation $\La^\nu \psi=\nu(2\lambda+n) \psi $, then $\lambda=l$ for
some $l=0,1,2, \cdots $, and we have
\[ [{{\Av}_{g_0}}\psi](z) = C e^{-\frac\nu 2|z|^2} {_1F_1}(-l; n;
\nu|z|^2)
\]
for some nonzero constant $C$ and certain $g_0\in G$.
    \end{itemize}
\end{proposition}

\noindent{\it Proof.} To prove i) let  $h\in U(n)$. Then
 \[
[{{\Av}_g}(\psi)](h.z)
 \stackrel{(\ref{average0})}
 {=} \int_{U(n)}
[T^\nu_{g k}(\psi)](h.z) dk = \int_{U(n)} j_\nu(g k,h.z)\psi(gkh.z)
dk
\]
Next, making use of the fact $j_\nu(g k,h.z)=j_\nu(g kh,z)$ for
$h,k\in U(n)$, the change $k'=kh\in U(n)$ and the $U(n)$-invariance
of the Haar measure $dk$ yield
\begin{align*} [{{\Av}_g}(\psi)](h.z)
 &=  \int_{U(n)} j_\nu(g
kh,z)\psi(gkh.z) dk\\
 &=  \int_{U(n)} j_\nu(g
k',z)\psi(gk'.z) dk \\ &=  [{{\Av}_g}(\psi)](z)
\end{align*}
 and therefore ${{\Av}_g}(\psi)$ is radial.

The assertion ii) on the boundedness of ${{\Av}_g}(\psi)$ follows
easily from the boundedness of the function $\psi$ keeping in mind
that $U(n)$ is compact.

The proof of iii) relies essentially to the fact that the Landau
Hamiltonian $\La^\nu$ commutes with the transformations $T^\nu_{g}$;
see Proposition \ref{Proposition31}.

For iv), let note first from the assumption that $\psi$ is a nonzero function on $\C^n$, there is $z_0\in \C^n$ such that
$\psi(z_0)\ne 0$ and then one can choose $g_0\in G$ such that
$g_0.0=z_0$ (such $g_0$ exists since the action of $G$ on $\C^n$ is
transitive). Hence using the fact that $k.0=0$ for every $k\in
U(n)$, it follows
\[[{{\Av}_{g_0}}(\psi)](0)= \psi(g_0.0)=\psi(z_0) \ne 0,\]
and therefore $ {{\Av}_{g_0}}(\psi)$ is a nonzero function on
$\C^n$. Now, the assertion in iv) follows by combining i), ii), iii)
above and i) and ii) of Proposition \ref{Proposition32}. \fin
\vspace*{.2cm}

 Thus, making use of the obtained properties of
the operators ${{\Av}_{g}}$ and $\K^{\nu}_{\Gamma,\chi}$, we show
the following

\begin{theorem}[Selberg transform]\label{SelbergTransformProp} 
Let $h$  be a given $L^{2,\nu}_{\Gamma,\chi}$-eigenvalue of the
integral operator $\K^{\nu}_{\Gamma,\chi}$. Then, there exists a
positive integer $j$ such that
\begin{equation}\label{SelbergTransform0}    h \varphi_j(z) =  \int_{\C^n}  e^{i\nu\omega(z,w)} Q^\nu(|z-w|)\varphi_j(w)
  dm(w), \end{equation}
where $\varphi_j(z):=e^{-\frac\nu 2|z|^2} {_1F_1}(-j; n; \nu
|z|^2).$ In particular
\begin{equation}\label{SelbergTransform}  h  =  \int_{\C^n}Q^\nu(|w|)
   {_1F_1}(-j; n; \nu
|w|^2) e^{-\frac\nu 2|w|^2}dm(w). \end{equation}
\end{theorem}

\noindent{\it Proof.} Since the operators $\K^{\nu}_{\Gamma,\chi}$
and $\La^\nu_{\Gamma,\chi}$ commute (see iii) of Proposition
\ref{KerFunCommu}), there is a basis of the Hilbert space
$L^{2,\nu}_{\Gamma,\chi}$ constituted of common eigenfunctions of
both $\K^{\nu}_{\Gamma,\chi}$ and $\La^\nu_{\Gamma,\chi}=\La^\nu$.
Hence, for $h$  being a given $L^{2,\nu}_{\Gamma,\chi}$-eigenvalue
of the integral operator $\K^{\nu}_{\Gamma,\chi}$, there exists
$\lambda\in\C$ and a nonzero function $\psi_\lambda \in
L^{2,\nu}_{\Gamma,\chi}$ such that
\begin{equation}
\label{Lem3Pr1} \K^{\nu}_{\Gamma,\chi}(\psi_\lambda ) =
h\psi_\lambda
\end{equation} and
\begin{equation}
\label{Lem3Pr2} \La^\nu(\psi_\lambda )= \nu(2\lambda+n) \psi_\lambda
.
\end{equation}
Hence, by taking the average of the involved functions in both sides
of (\ref{Lem3Pr1}) and using
 Proposition \ref{PropSelberg1}, we
deduce
\begin{align*}
  h  [{{\Av}_{g_0}} (\psi_\lambda)](z) &=  {{\Av}_{g_0}} [\K^{\nu}_{\Gamma,\chi}
  (\psi_\lambda)](z) =
  {{\Av}_{g_0}} [\K
  (\psi_\lambda)](z) =
  {\K}[{{\Av}_{g_0}}(\psi_\lambda)](z)
\end{align*}
where ${g_0}$ is certain element of $G$ such that
$\psi_\lambda({g_0}.0)\ne 0$. Hence, we have
\begin{equation*}
  h  [{{\Av}_{g_0}} (\psi_\lambda)](z)
  = \int_{\C^n} e^{i\nu\omega(z,w)}Q^\nu(|z-w|)
  [{{\Av}_{g_0}} (\psi_\lambda)](w)dm(w).
\end{equation*}
 In the other hand it
follows from iv) of Proposition \ref{AveragingFunction} that
$\lambda =j$ for certain positive integer $j$ and that
${{\Av}_{g_0}} (\psi_\lambda)$ is proportional to the radial
 solution $\psi_j$ of the
differential equation $\La^\nu \phi=\nu(2j+n) \phi$, i.e.,
\[ [{{\Av}_{g_0}} (\psi_j)](w) = C\varphi_j(w) := C e^{-\frac\nu 2|w|^2} {_1F_1}(-j; n; \nu
|w|^2)\] with  $C= [{{\Av}_{g_0}} (\psi_j)](0)=
 \psi_j(g_0.0) \ne 0$. Therefore, we obtain
\[
  h \varphi_j(z) =  \int_{\C^n}  e^{i\nu\omega(z,w)} Q^\nu(|z-w|)\varphi_j(w)
  dm(w).
\]
 In particular, we get the expression of $h$ as
 asserted in (\ref{SelbergTransform}) when taking $z=0$.
This completes the proof of Theorem \ref{SelbergTransformProp}. \fin

\section{Proof of main result. \label{Proofs} }

First let recall that, under the $(RDQ)$ condition, the space
${\mathcal{F}}^\nu_{\Gamma,\chi}$ can be identified to the space of
$\mathcal{C}^\infty$ sections of an appropriate line bundle
$\mathcal{L}_{\nu,\gamma,\chi}$ over the compact manifold
$\C^n/\Gamma\cong \Lambda(\Gamma)$ (see Remark 3.6).
  Hence, it is a standard fact \cite{Gilkey} that the self-adjoint elliptic
differential operator $\La^\nu_{\Gamma,\chi}$ has discrete spectrum
$Sp(\La^\nu_{\Gamma,\chi})$, constituted of an increasing sequence
of eigenvalues tending to $+\infty$ and occuring with finite
degeneracy. Furthermore, we have the following orthogonal
decomposition
\begin{equation}\label{OrthDecomp} L^{2,\nu}_{\Gamma,\chi}=\bigoplus_{\mu(\lambda)\in
Sp(\La^\nu_{\Gamma,\chi})} {\mathcal{E}}^\nu_{\Gamma,\chi}(\lambda)
\end{equation}
 with $\mu(\lambda)=\nu(2\lambda +n)$ and where
 ${\mathcal{E}}^\nu_{\Gamma,\chi}(\lambda)$; $\lambda\in\C$, are the eigenspaces
\[{\mathcal{E}}^\nu_{\Gamma,\chi}(\lambda)=\set{ f\in
{\mathcal{F}}^\nu_{\Gamma,\chi};
  \,\, \La^\nu_{\Gamma,\chi}f=\nu(2\lambda+n)f}.\]
  Thus the main result of this paper concerns the concrete
description of the orthogonal decomposition (\ref{OrthDecomp}).
Mainly, we determine the spectrum of the Landau Laplacian
$\La^\nu_{\Gamma,\chi}$ when acting on $L^{2,\nu}_{\Gamma,\chi}$ and
we show that the eigenspaces
${\mathcal{E}}^\nu_{\Gamma,\chi}(\lambda)$ corresponding to $\lambda
=l=0,1,2, \cdots, $ are the only nonzero (finite) eigenspaces,
whose the dimension is computed explicitly. Also, we give the
explicit formula for their reproducing kernels. To do this, we begin
by
 fixing a given positive integer $l$ and considering the eigenprojector kernel of the
$L^2$-eigenspace $\mathcal{A}^{2,\nu}_{l}(\C^n)$ as given in
Proposition \ref{Proposition32} by (\ref{RepK}), to wit
\[{\mathcal{K}}^\nu_l(z,w) = e^{i\nu\omega(z,w)}Q^\nu_l(|z-w|) ,\]
    with
\[ Q^\nu_l(|z-w|) =\frac{\Gamma(n+l)}{\Gamma(n)l!}
({\nu}/{\pi})^ne^{-\frac\nu 2|z-w|^2} {_1F_1}(-l; n; \nu |z
-w|^2).\] Let ${\mathcal{K}}^\nu_{l;\Gamma,\chi}(z,w)$ be its
associated $(\Gamma,\chi)$-automorphic kernel function. Then from
the previous section, we know that there is  $w_0 \in \C^n$ such
that $ z \longmapsto {\mathcal{K}}^\nu_{l;\Gamma,\chi}(z,w_0)$ is a
nonzero function on $\C^n$ belonging to the space
${\mathcal{F}}^\nu_{\Gamma,\chi}$ and furthermore to the eigenspace
${\mathcal{E}}^\nu_{\Gamma,\chi}(l)$. Thus, we assert

\begin{proposition} \label{Prop52} Assume the $(RDQ)$ condition to be satisfied by the triplet
$(\nu,\Gamma, \chi)$. Then the eigenspace
${\mathcal{E}}^\nu_{\Gamma,\chi}(\lambda)$ is a nonzero vector
space if and only if  $\lambda=l$ is a positive integer $l=0,1,2,
\cdots$.
\end{proposition}

\noindent{\it Proof.} We have to prove only the "only if". Indeed,
it holds by applying iv) of Proposition \ref{AveragingFunction} to a
given nonzero function
$\psi\in{\mathcal{E}}^\nu_{\Gamma,\chi}(\lambda)\ne \set{0}$.
 \fin\\

The dimension formula of ${\mathcal{E}}^\nu_{\Gamma,\chi}(l)$,
$l=0,1,2, \cdots,$ is an immediate consequence of the following

\begin{proposition}\label{Prop51}
Fix $l=0,1,2, \cdots$ and let $\K^\nu_{l;\Gamma,\chi}$ be the
integral operator on $L^{2,\nu}_{\Gamma,\chi}$ associated with the
$(\Gamma,\chi)$-automorphic kernel function
${\mathcal{K}}^\nu_{l;\Gamma,\chi}(z,w)$,
 \begin{equation}\Big[\K^\nu_{l;\Gamma,\chi}(\psi)\Big](z) =
\int_{\C^n/\gamma} {\mathcal{K}}^\nu_{l;\Gamma,\chi}(z,w) \psi(w)
dm(w) .\end{equation}
 Then the eigenvalues $h_{l,j}$;
$j=0,1,2, \cdots $ of $\K^\nu_{l;\Gamma,\chi}$ are given explicitly
by
\[h_{l,j} = \left \{\begin{array}{ll} 1 & \quad \mbox{if} \quad j=l \\ 0 & \quad \mbox{otherwise} \end{array}\right. .\]
\end{proposition}

 \noindent{\it Proof.}
According to Theorem \ref{SelbergTransformProp}, it follows
 that if $h_{l,j}$ is an $L^{2,\nu}_{\Gamma,\chi}$-eigenvalue  of
 $\K^\nu_{l;\Gamma,\chi}$, then it is given through (\ref{SelbergTransform}) by
  \begin{align*} h_{l,j}  &= \int_{\C^n}Q^\nu_l(|w|)
   {_1F_1}(-k; n; \nu
|w|^2) e^{-\frac\nu 2|w|^2}dm(w)\\ &=
\frac{\Gamma(n+l)}{\Gamma(n)l!}({\nu}/{\pi})^n
\int_{\C^n}{_1F_1}(-l; n; \nu |w|^2)
   {_1F_1}(-k; n; \nu
|w|^2) e^{-\nu |w|^2}dm(w)\end{align*} for  some (unique) positive
integer $k=k_j$.
 Next, by the use of the polar
coordinates $z=r\theta$, $r\geq 0$, $\theta\in S^{2n-1}$, the change
of variable $x=\nu r^2$ and the fact that the involved
hypergeometric function ${_1F_1}(-j;c;x)$
 for $j\in \Z^+$ is related to the orthogonal Laguerre polynomial
 $L^{c-1}_j(x)$ \cite[page 333]{NO},
  \[  {_1F_1}(-j;c;x) = \frac{j!\Gamma(c)}{\Gamma(j+c)}
  L^{c-1}_j(x),\]
   we get
    \[ h_{l,j} =\frac{vol(S^{2n-1})}{2 \pi^n}
\frac{k!\Gamma(n)}{\Gamma(k+n)}
     \int_0^{+\infty} L^{n-1}_l(x)L^{n-1}_{k}(x)x^{n-1}e^{-x}
dx.\] The above involved integral is known to be given by \cite[page
56]{NO}, \[ \int_0^{+\infty} L^{n-1}_l(x)L^{n-1}_{k}(x)x^{n-1}e^{-x}
dx=\frac{\Gamma(l+n)}{l!}\delta_{kl}.\] Thus, since
$\mbox{vol}(S^{2n-1}) = 2\pi^{n}/\Gamma(n)$, we get $h_{l,j} =
\delta_{kl},$ where $\delta_{kl}$ denotes here the Kronecker symbol.
\fin\vspace*{.2cm}

Using the above result together with (\ref{SelbergTransform0}) in
Theorem \ref{SelbergTransformProp}, keeping in mind that
$Q^\nu_l(0)=1$, one deduces easily the following
\begin{corollary}
The function ${\mathcal{K}}^\nu_{l;\Gamma,\chi}(z,w)$ is the
reproducing kernel of the eigenspace
${\mathcal{E}}^\nu_{\Gamma,\chi}(l)$, that is $$ f(z) =
\int_{\C^n/\Gamma}{\mathcal{K}}^\nu_{l;\Gamma,\chi}(z,w) f(w)dm(w)$$
for every $f\in {\mathcal{E}}^\nu_{\Gamma,\chi}(l)$. Therefore, we
have
\begin{align*} \dim
{\mathcal{E}}^\nu_{\Gamma,\chi}(l)=\int_{\C^n/\Gamma}{\mathcal{K}}^\nu_{l;\Gamma,\chi}(z,z)
dm(z)=\frac{\Gamma(n+l)}{\Gamma(n)l!}({\nu}/{\pi})^n\mbox{vol}(\Lambda(\Gamma)).
\end{align*}
\end{corollary}

Hence to complete the proof of the main theorem, we have only to
prove the following
\begin{proposition} \label{Prop53} Under the $(RDQ)$ condition, the functional
 space ${\mathcal{O}}^\nu_{\Gamma,\chi}(\C^n)$
defined by (\ref{HolAutForms}) is isomorphic to the fundamental
space of $(\Gamma,\chi)$-ground states,
 ${\mathcal{E}}^\nu_{\Gamma,\chi}(0)$. Precisely,
 $f \longmapsto g=e^{\frac\nu 2|z|^2}f$ defines an isomorphism map from
${\mathcal{E}}^\nu_{\Gamma,\chi}(0)$ onto
${\mathcal{O}}^\nu_{\Gamma,\chi}(\C^n)$.
\end{proposition}

\noindent{\it Proof.}
 The assertion is a
consequence of Proposition \ref{Prop51} combined with iii) of
Corollary \ref{LastRepKer}.
  Indeed, from the
explicit obtained dimensional formulas, we deduce that  $\dim
{\mathcal{O}}^\nu_{\Gamma,\chi}(\C^n) =\dim
{\mathcal{E}}^\nu_{\Gamma,\chi}(0).$ Therefore, making the
observation that $$
\mathfrak{G}^{-1}[{\mathcal{O}}^\nu_{\Gamma,\chi}(\C^n)] \subset
{\mathcal{E}}^\nu_{\Gamma,\chi}(0)$$ for
 ${\mathcal{O}}^\nu_{\Gamma,\chi}(\C^n)\subset \ker \Delta^\nu_{\Gamma,\chi}
 = \mathfrak{G}[{\mathcal{E}}^\nu_{\Gamma,\chi}(0)],$
it follows also that
\[{\mathcal{O}}^\nu_{\Gamma,\chi}(\C^n) = \ker \Delta^\nu_{\Gamma,\chi} \cong
  {\mathcal{E}}^\nu_{\Gamma,\chi}(0).\]
This completes the proof. \fin

\begin{remark}
According to the proof of the last proposition, we note that the
space of holomorphic functions on $\C^n$ satisfying the functional
equation $$g(z+\gamma) = \chi(\gamma) e^{\frac \nu 2
|\gamma|^2+\nu\scal{z,\gamma}}g(z)$$ can be viewed as the null space
of the differential operator
$$\Delta^\nu_{\Gamma,\chi}=\sum\limits_{j=1}\limits^n(
\frac{-\partial^2}{\partial z_j\partial \bar z_j} + \nu \bar z_j
\frac{\partial}{\partial \bar z_j})$$ acting on  $\mathcal C^\infty$
functions on ${\mathcal{G}}^\nu_{\Gamma,\chi}$. This means that the
number of first order differential operators defining such
holomorphic functions reduces further to the single elliptic second
order differential operator $\Delta^\nu_{\Gamma,\chi}$. Such
situation arise frequently in the theory of holomorphic functions.
\end{remark}

\section{Concluding remarks}

Using the notation of Section 3 (the Table there) and assuming the
$(RDQ)$ condition to be satisfied, we see that the main theorem can
be reworded in terms of the spectral properties of the differential
operator
\[ \Delta^\nu_{\Gamma,\chi}= \sum_{j=1}^n \Big(
\frac{-\partial^2}{\partial z_j \partial \bar z_j}+ \nu
 \bar z_j \frac{\partial}{\partial
\bar z_j} \Big)\] acting on the Hilbert space
$B^{2,\nu}_{\Gamma,\chi}=\overline{{\mathcal{G}}^\nu_{\Gamma,\chi}}
     ^{\scal{\scal{\cdot , \cdot }}_\Gamma}$ endowed with the
Hermitian inner product (\ref{HSPGamma2}),
\begin{align}
\scal{\scal{g_1, g_2}}_\Gamma=\int_{\C^n/\Gamma}g_1(z)
\overline{g_2}(z) e^{-\nu|z|^2}dm(z).\label{InnerSP}\end{align}
 Namely, we have the
following Hilbertian orthogonal decomposition
\[B^{2,\nu}_{\Gamma,\chi}= \bigoplus_{l=0}^\infty
E_l(\Delta^\nu_{\Gamma,\chi}), \] where
$E_l(\Delta^\nu_{\Gamma,\chi})$, $l=0,1, \cdots ,$ are the
eigespaces defined by $$E_l(\Delta^\nu_{\Gamma,\chi})=\set{g\in
{\mathcal{G}}^\nu_{\Gamma,\chi}; \quad  \Delta^\nu_{\Gamma,\chi} g =
\nu l g}.$$

 According to the proof of ii) of the main theorem, the eigenspace
$E_0(\Delta^\nu_{\Gamma,\chi})$, corresponding to $l=0$, is nothing
but the space ${\mathcal{O}}^\nu_{\Gamma,\chi}(\C^n)$ studied in
Section 3. Thus, the space of holomorphic functions on $\C^n$,
$n\geq 1$, satisfying (\ref{ChiGammaPeriodicEntire}), can then be
realized as the null space of the differential operator
$\Delta^\nu_{\Gamma,\chi}$ on ${\mathcal{G}}^\nu_{\Gamma,\chi}$.
Also, let note here that for $\nu=\pi$,
${\mathcal{O}}^\pi_{\Gamma,\chi}(\C^n)$ is called the space of
$\Gamma$-Theta functions with respect to $(H,\chi)$;
$H(z,w)=\scal{z,w}$, and plays an important role in the theory of
abelian varieties \cite{Mumford}.

Now, since ${\mathcal{O}}^\nu_{\Gamma,\chi}(\C^n)=Ker
(\Delta^\nu_{\Gamma,\chi})$ and its dimension is given by
$$\dim{\mathcal{O}}^\nu_{\Gamma,\chi}(\C^n)=({\nu}/{\pi})^n\mbox{vol}(\C^n/\Gamma),$$ we may
conclude this paper by asking if there is a some how canonical way
for construction of an {\it orthogonal basis} of
${\mathcal{O}}^\nu_{\Gamma,\chi}(\C^n)$ with respect to the
Hermitian inner product (\ref{InnerSP}). We may address the same
question for the other eigenspaces $E_l(\Delta^\nu_{\Gamma,\chi})$,
$l=0,1,2, \cdots,$  whose the dimension is given by
$$\dim E_l(\Delta^\nu_{\Gamma,\chi})=\frac{\Gamma(n+l)}{\Gamma(n)l!}({\nu}/{\pi})^n\mbox{vol}(\C^n/\Gamma).$$
For the case $n=1$ and $l=0$, Professor Y. Hantout
\footnote{University of Lille 1, France.} has communicated to us the
construction of
 a canonical orthogonal basis of $({\mathcal{O}}^\nu_{\Gamma,\chi}(\C),
\scal{\scal{\cdot , \cdot}}_\Gamma)$ expressed in terms of some
Jacobi Theta functions in one variable $z$. We hope to handle the
general case ${\mathcal{O}}^\nu_{\Gamma,\chi}(\C^n;H)$ in a near
future.\\

{\bf\it Acknowledgements.} The authors would like to thank the organizers
 of the meeting AHGS06 ({\it Moroccan Association of Harmonic Analysis and
 Spectral Geometry}), where a part of this work was presented as talk.
 A. G. would like to address special thanks to the Center for Advanced
 Mathematical Sciences (CAMS) of the American University of
 Beirut (AUB) for the hospitality during the full academic year 2006-2007.

\end{document}